\documentclass[12pt]{article}
\usepackage[utf8]{inputenc}
\usepackage[T2A]{fontenc}
\usepackage[utf8]{inputenc}
\usepackage{amsthm,amssymb,amsmath}
\usepackage{dsfont}
\usepackage{t1enc}
\usepackage[english]{babel}
\usepackage{xcolor}
\usepackage{cmap}
\usepackage{mathrsfs}

\numberwithin{equation}{section}

\newtheorem{con}{Convention}

\newtheorem{theorem}{Theorem}[section]
\newtheorem{lemma}[theorem]{Lemma}
\newtheorem{rem}[theorem]{Remark}
\newtheorem{cor}[theorem]{Corollary}
\newtheorem{assump}[theorem]{Assumption}
\newtheorem{defin}[theorem]{Definition}
\newtheorem{prop}[theorem]{Proposition}

\newtheorem*{thm*}{Theorem}

\begin{document}

\bigskip

\title{Anomalous singularity of the solution of the vector Dyson equation in the critical case\footnote{This work was done under the summer internship 
program of IST Austria. }}
\author{Oleksii Kolupaiev\footnote{avkolupaev0509@gmail.com},\\ 
V.N. Karazin Kharkiv National University}
\date{\today}
\maketitle

{\footnotesize
\begin{center}
\textbf{Abstract}\\
\end{center}
We consider the solution of the vector Dyson equation $-1/m=z+Sm$ in the case when $S$ has a block staircase structure with $(n-1)$ different critical zero blocks below the strictly positive anti-diagonal and all elements right above the anti-diagonal are strictly positive. We prove that the components of $m$ behave as fractional powers of $z$ in the neighbourhood of zero and show that the self-consistent density of states $\rho(E)$ behaves as $\vert E\vert^{-\frac{n-1}{n+1}}$ as $E$ tends to zero, where $n^2$ is a number of blocks. Both constant block and non-constant block cases are considered. In the non-constant case uniform estimates for the components of $m$ are obtained.\\  

}



\section{Set-up and main results.}
Consider a generalized Wigner matrix $H$, i.e. a hermitian random matrix $H=\left(h_{i,j}\right)_{i,j=1}^n$ such that $\mathbb{E} H=0$, $h_{i,j}$ are independent up to a hermitian symmetry, but not necessarily identically distributed. 
One way to describe the distribution of eigenvalues of $H$ is to calculate the trace of the resolvent $G=\left(H-z\right)^{-1}$ and to apply the inverse Stieltjes transform to it (e.g. see [1]). The diagonal entries of $G$ can be typically well approximated by the components of the solution of the vector Dyson equation (VDE)\\  
\begin{equation}
-\frac{1}{m}=z+Sm,
\label{VDE}
\end{equation}\\
these results are called entrywise (local) laws. In \eqref{VDE} $S=\left(\mathbb{E}\left\vert h_{i,j}\right\vert^2\right)_{i,j=1}^n$ is a real-symmetric matrix with non-negative entries;  $z \in \mathbb{H}$ (we denote complex upper half-plane by $\mathbb{H}$); $m$ is a vector-function, $m: \mathbb{H}\rightarrow \mathbb{H}^n$. Here $-1/m$ is understood in entry-wise sense, by $z$ we denote a vector from $\mathbb{H}^n$ whose components are equal to $z$. From [1,~Theorem~6.1.4] it is known that (\ref{VDE}) has a unique solution. The self-consistent density of states is defined in [1]:\\
\begin{equation}
\rho(E)=\lim_{\eta\searrow 0}{\frac{1}{\pi n}\sum\limits_{k=1}^n \operatorname{Im}{m_k\left(E+i\eta\right)}}. 
\end{equation}

We impose the following assumption on the matrix $S$: $S$ does not have any rectangular zero blocks with a perimeter at least $2(n+1)$. It is known that in the case when $S$ does not meet this assumption, $\|m(z)\|$ behaves as $\frac{1}{\vert z\vert}$ near zero and a delta-measure appears in $\rho(E)$ at zero. This reflects the fact that the corresponding random matrix $H$ does not have full rank. On the other hand, if $S$ has only zero blocks with a perimeter strictly less than $2n$ (the situation extensively investigated in [4]), then the solution of \eqref{VDE} is bounded. We will consider the case, when $S$ has a staircase structure with $(n-1)$ different critical zero block (i.e. with a perimeter exactly $2(n+1)$) below the strictly positive anti-diagonal. We also assume strict positivity of the elements right above the anti-diagonal, namely\\

\begin{assump}\label{assump_S}
For all $i,j\in\left\lbrace 1,\ldots,n\right\rbrace$ we have\\
\begin{equation*}
\begin{cases}
s_{i,j}>0,\text{if}\quad i+j=n\\
s_{i,j}>0,\text{if}\quad i+j=n+1\\
s_{i,j}=0,\text{if}\quad i+j\ge n+2\\
\end{cases}
\end{equation*}
\end{assump}

We also adhere to the following conventions:\\
\begin{con} The variable $z$ always lies in the complex upper half-plane. Particularly, if we consider $z\rightarrow 0$ it means that $z$ tends to zero from the upper half-plane.
\end{con}

\begin{con} $c,C$ are positive constants, which can change their value from line to line for simplicity of notation.
\end{con}

\begin{defin}\label{def_equiv}
Let $\left\lbrace z_l\right\rbrace_{l=1}^\infty\subset\mathbb{H}$ be a fixed sequence converging to zero, $f,g:\mathbb{H}\rightarrow \mathbb{C}$. We will say that $f$ is equivalent to $g$ with respect to the sequence $\left\lbrace z_l\right\rbrace_{l=1}^\infty$ if there exist constants $c,C$ such that $c<\frac{\vert f(z_l)\vert}{\vert g(z_l)\vert}<C$ for all $l\in\mathbb{N}$.
Equivalence of functions $f$ and $g$ will be denoted by $f\sim g$. 
\end{defin}

Our main results are as follows:\\

\begin{theorem}\label{t_main}
Under Assumption \ref{assump_S} there exist positive constants $\left\lbrace c_k\right\rbrace_{k=1}^n$ such that for all $k\in\left\lbrace 1,\ldots,n\right\rbrace$ we have\\
\begin{equation}
\lim_{z\to 0} {\left(m_k(z)z^{-\left(1-\frac{2k}{n+1}\right)}\right)} = c_ke^{i\frac{\pi k}{n+1}},
\label{f_main}
\end{equation}
where the branch of $z^{-\left(1-\frac{2k}{n+1}\right)}$ is chosen in the complex upper half-plane in such a way that $z^{-\left(1-\frac{2k}{n+1}\right)}$ is positive for positive $z$.\\ 
\end{theorem}

\begin{cor}\label{cor_main}
Under assumptions and with notations of Theorem \ref{t_main} we have that the following limit exists and it is positive\\
\begin{equation}
\lim_{E\to 0} {\rho(E)\vert E\vert ^{-\frac{n-1}{n+1}}}>0.
\end{equation}
\end{cor}

A weaker version of Corollary \ref{cor_main} was already obtained in [3]: the matrix $S$ was considered under the stronger assumption $s_{i,j}>0$ for all $i+j\le n+1$ and only the equivalence $\rho (E)\sim \vert E\vert ^{-\frac{n-1}{n+1}}$, as $E$ tends to zero, was established without proving the limit.\\

The equation \eqref{VDE} typically arises when $H$ is a random matrix consisting of $n\times n$ blocks of size $N\times N$ each. We think of $n$ to be fixed and $N\to\infty$. If the variances inside each block are constant $s_{i,j}$, $i,j\in\left\lbrace 1,\ldots,n\right\rbrace$, then $G_{l,l}(z)\approx m_k(z)$ if $l$ belongs to the $k$-th block and $m_k$ is the $k$-th component of the solution \eqref{VDE}. One may also consider the case of non-constant blocks in the matrix of variances $S$, where the corresponding VDE has $nN$ possible different components. Nevertheless, the following result similar to Theorem \ref{t_main} still holds:\\

\begin{theorem}\label{t_main_non_const}\emph{(Non-constant block case.)}\\ 
Let $\mathcal{S}$ be a real-symmetric $nN\times nN$ matrix with non-negative entries. Assume that $\mathcal{S}$ has a following block structure:\\
\begin{equation}
\mathcal{S}=\left(\mathcal{S}^{j,k}\right)_{j,k=1}^n,
\end{equation}
where $\mathcal{S}^{j,k}$ are $N\times N$ matrices for $j,k\in\left\lbrace 1,\ldots,n\right\rbrace$ and\\
\begin{enumerate}
\item{ All entries of $\mathcal{S}^{j,k}$ are either positive or all entries are zero.}\\
\item{ If all entries of $\mathcal{S}^{j,k}$ are positive then $\frac{\mathfrak{c}}{N}<\mathcal{S}^{j,k}_{\nu,\tau}<\frac{\mathfrak{C}}{N}$ for all $\nu,\tau\in\left\lbrace 1,\ldots,N\right\rbrace$, where $\mathfrak{c}$ and $\mathfrak{C}$ are positive constants which do not depend on $N$.}
\item{ $\mathcal{S}^{j,k}=0$ if $j+k>n+1$.}\\
\item{ All entries of $\mathcal{S}^{j,k}$ are positive if $j+k\in\left\lbrace n,n+1\right\rbrace$.}\\
\end{enumerate}
Then there exist positive constants $\mathscr{C}_1$,$\mathscr{C}_2$ which do not depend on $N$ and a neighbourhood $U$ of zero in $\mathbb{H}$ such that for all $k\in\left\lbrace 1,\ldots,nN\right\rbrace$ and $z\in U$ it holds that\\
\begin{equation}
\mathscr{C}_1< \vert m_k(z)\vert \cdot \vert z\vert ^{-\left(1-\frac{2}{n+1}\left\lceil\frac{k}{N}\right\rceil\right)}< \mathscr{C}_2.
\label{f_main_non_const_abs}
\end{equation} 
Moreover,\\

\begin{equation}
\operatorname{arg}{\left(m_k(z)z^{-\left(1-\frac{2}{n+1}\left\lceil\frac{k}{N}\right\rceil\right)}\right)} \to \frac{\pi}{n+1}\left\lceil\frac{k}{N}\right\rceil, \quad \text{as}\quad z\to 0, z\in\mathbb{H},
\label{f_main_non_const_ph}
\end{equation}
where the branch of $z^{-\left(1-\frac{2}{n+1}\left\lceil\frac{k}{N}\right\rceil\right )}$ is chosen in the complex upper half-plane in such a way that $z^{-\left(1-\frac{2}{n+1}\left\lceil\frac{k}{N}\right\rceil\right)}$ is positive for positive $z$.\\ 

\end{theorem}

Note that the arguments of the solution still converge within each block \eqref{f_main_non_const_ph}, but their moduli typically do not converge. However we have uniform bounds on their asymptotics \eqref{f_main_non_const_abs}.\\

In Sections 2 to 6 we prove Theorem \ref{t_main} and in Section 7 prove Theorem \ref{t_main_non_const}  generalizing arguments from the previous sections to the non-constant block case.\\
We prove Theorem \ref{t_main} in two steps: (\ref{f_main}) is verified separately for absolute value and for argument of the function $m_k(z)z^{-\left(1-\frac{2k}{n+1}\right)}$.\\

During the preparation of this article D.Renfrew informed us on an independent ongoing work [5] studying similar questions for the same model. In [5] the classification of singularities of the solution at zero is given for any constant block  matrix S in contrast to our Assumption 1.1 but the non constant block case (Theorem 1.5) was not considered in [5]. Note that the method presented in [5] is completely different from the method described in the current paper.\\

\textbf{Acknowledgement:} This work was done during the summer internship at L\'aszl\'o Erd{\H o}s' research group at IST Austria. I am very grateful to him for supervising this project. Also I would like to thank Asbjorn Lauristen for kindly discussing [3].\\

\section{Relation between the behaviour of the solution of (\ref{VDE}) and the structure of $S$.}\label{sec2}

We fix a sequence $\left\lbrace z_l\right\rbrace_{l=1}^{\infty}\subset \mathbb{H}$ such that\\
\begin{equation}
\lim_{l\to\infty}z_l =0.
\end{equation} 
In Sections 2 to 4 we will temporarily need the following technical assumption, which will be removed in Section 5:\\

\begin{assump}\label{assump_zm}
\begin{equation}
\lim_{l\to\infty} {\left(z_l\cdot m_j(z_l)\right)} = 0,\quad \forall j\in\left\lbrace 1,\ldots,n\right\rbrace.
\end{equation}
\end{assump}

In Subsection \ref{subs_staircase} we establish a block staircase structure of the matrix $S$ assuming that components of the solution of (\ref{VDE}) are ordered in the sense of their asymptotic size near 0. First we formulate this ordering as an assumption and later in Subsection \ref{subs_sol_ord} we show that it follows from Assumption \ref{assump_S}.\\

\begin{assump}\label{assump_sol_ord}
There exists a subsequence $\left\lbrace z_{l_q}\right\rbrace_{q=1}^\infty$ of the sequence $\left\lbrace z_l\right\rbrace_{l=1}^\infty$ such that for every $k\in\left\lbrace 1,\ldots,n-1\right\rbrace$ either $m_k\sim m_{k+1}$ with respect to the sequence $\left\lbrace z_{l_q}\right\rbrace_{q=1}^{\infty}$ or\\
\begin{equation}
\lim_{q\to\infty} {\frac {m_k(z_{l_q})}{m_{k+1}(z_{l_q})}} = 0.
\end{equation}
\end{assump}

\subsection{Block staircase structure of the matrix $S$.}\label{subs_staircase}
In this subsection we work under Assumptions \ref{assump_zm} and \ref{assump_sol_ord}. For simplicity of notation we denote $w_q:=z_{l_q}$, $q\in\mathbb{N}$.\\

Let us split the vector $m$ into $p$ blocks $\left(m^{(1)},m^{(2)},\ldots ,m^{(p)}\right)$ of sizes $d_1,d_2,\ldots,d_p$, $m^{(a)}\in\mathbb{C}^{d_a}$ for $a\in\left\lbrace 1,\ldots,p\right\rbrace$, in such a way that two components of $m$ belong to the same block if and only if they are equivalent with respect to the sequence $\left\lbrace w_q\right\rbrace_{q=1}^\infty$, i.e. for $j,k\in\left\lbrace 1,\ldots,n\right\rbrace$
\begin{equation*}
m_j\sim m_k \quad \text{WRT} \quad \left\lbrace w_q\right\rbrace_{q=1}^\infty \Leftrightarrow 
\end{equation*}
\begin{equation}
\Leftrightarrow\exists a\in\left\lbrace 1,\ldots,p\right\rbrace; \exists \nu,\tau\in\left\lbrace 1,\ldots,d_a\right\rbrace:\quad j = \sum\limits_{b=1}^{a-1}d_b + \nu,\quad k = \sum\limits_{b=1}^{a-1}d_b+\tau.
\end{equation}
In short, all components within each block are comparable along $\left\lbrace w_l\right\rbrace_{l=1}^\infty$, and\\
\begin{equation}
\|m^{(a)}\|\ll \| m^{(b)}\| \quad \text{iff} \quad a<b
\label{f_diff_blocks}
\end{equation}
by Assumption \ref{assump_sol_ord}.\\

Now a corresponding partition of the matrix $S$ into a block matrices is naturally obtained: $S=\left(\mathfrak{S}_{a,b}\right)_{a,b=1}^p$, where $\mathfrak{S}_{a,b}$ is a rectangular matrix of the size $d_a\times d_b$. As $S$ is symmetric, $\mathfrak{S}_{a,b}^T=\mathfrak{S}_{b,a}$, $a,b=1,\ldots ,p$. In the notation above the initial VDE (\ref{VDE}) takes the form:\\
\begin{equation}
-\frac {1}{m^{(a)}} = z + \sum\limits_{b=1}^p \mathfrak{S}_{a,b}m^{(b)}, \quad a=1,\ldots ,p,
\label{nfVDE}
\end{equation}
or, equivalently,\\
\begin{equation}
-1 = zm^{(a)} + \sum\limits_{b=1}^p m^{(a)}\mathfrak{S}_{a,b}m^{(b)}, \quad a=1,\ldots ,p.
\label{nfVDE1}
\end{equation} 

We will prove that\\
 
\begin{prop}\label{prop_staircase}
Under Assumption \ref{assump_sol_ord} the matrix $S=\left(\mathfrak{S}_{a,b}\right)_{a,b=1}^p$ has a block staircase structure in the following sense:\\
\begin{enumerate}
\item{$\mathfrak{S}_{a,p-a+1}$ is a non-zero matrix for $a\in\left\lbrace 1,\ldots, p\right\rbrace$,}\\
\item{$\mathfrak{S}_{a,b}=0$ for all $a,b\in\left\lbrace 1,\ldots,p\right\rbrace$ such that $a+b>p+1$.}
\end{enumerate}
Additionally under Assumption \ref{assump_zm} we also have that $\mathfrak{S}_{a,p-a+1}$ is a square matrix for $a\in\left\lbrace 1,\ldots, p\right\rbrace$.\\
\end{prop}

Define $r_a$ to be the index of the rightmost not identically zero block in the $a$-th row of blocks, i.e.:\\
\begin{equation}
r_a:=\max\left\lbrace b\mid \mathfrak{S}_{a,b}\neq 0\right\rbrace, \quad a=1,\ldots ,p.
\end{equation}

\begin{lemma} \label{l_summand_order_1}
Under Assumption \ref{assump_sol_ord} for every $a,b\in\left\lbrace 1,\ldots, p\right\rbrace$ and $\nu\in\left\lbrace 1,\ldots, d_a\right\rbrace$, $\tau\in\left\lbrace 1,\ldots, d_b\right\rbrace$ it holds that:\\
\begin{enumerate}
\item{$m^{(a)}_{\nu}(w_q)m^{(b)}_{\tau}(w_q)\to 0$ as $q\to\infty$ if $b<r_a$,}\\
\item{$m^{(a)}_{\nu}m^{(b)}_{\tau}\sim 1$ WRT $\left\lbrace w_q\right\rbrace_{q=1}^\infty$ if $b=r_a$.}\\
\end{enumerate}
\end{lemma}

In order to prove Lemma \ref{l_summand_order_1} we will need the notion of saturated self-energy operator from [1], that is a linear operator defined by the matrix\\
\begin{equation}
F=\vert m\vert S\vert m\vert = \left(\vert m_k \vert s_{k,j} \vert m_j\vert\right)_{k,j=1}^{n}.
\end{equation}
Further we will denote for short $\frac{m}{\vert m\vert} := \left\lbrace \frac{m_j}{\vert m_j\vert}\right\rbrace _{j=1}^n$.\\

From [1, Proposition 7.2.9] it is known that\\
\begin{equation}
\| F\|_{L^2\rightarrow L^2} < 1,
\label{f_F_estimate}
\end{equation} 
where $L^2 = \left(\mathbb{C}^n, \|.\|\right)$, for $v=\left(v_1,\ldots,v_n\right)\in \mathbb{C}^n$:\\ 
\begin{equation*}
\|v\|=\left(\frac{1}{n}\sum\limits_{j=1}^n\vert v_j\vert^2\right)^{\frac{1}{2}}.
\end{equation*}

From \eqref{f_F_estimate} it follows that for all $j,k\in\left\lbrace 1,\ldots, n\right\rbrace$ and $z\in\mathbb{H}$
\begin{equation}
\vert m_k(z)\vert s_{k,j}\vert m_j(z)\vert < 1.
\label{f_F_entrywise}
\end{equation}

We will also need the following estimates for the solution of \eqref{VDE}.\\
 
\begin{prop}\label{prop_sol_bound}
For all $k=1,\ldots,n$ and $z\in\mathbb{H}$, $\vert z\vert <1$ we have:\\
\begin{equation}
c\vert z\vert < \vert m_k(z)\vert < \frac{C}{\vert z\vert}.
\label{f_sol_bound}
\end{equation}
\end{prop}

\textbf{Proof of Proposition \ref{prop_sol_bound}:} From [1, Proposition 7.2.9] it is known that
\begin{equation}
\|m(z)\|_2\le\frac{2}{\vert z\vert}, \quad \forall z\in\mathbb{H}.
\end{equation}
Therefore, the upper estimate in \eqref{f_sol_bound} holds. On the other hand,\\
\begin{multline*}
-\frac{1}{m(z)}=z+Sm(z) \Rightarrow \left\|\frac{1}{m(z)}\right\|_2 = \|z+Sm(z)\|_2\le \|z\|_2 +\|S\|\|m(z)\|_2\le \vert z\vert +\frac{2\|S\|}{\vert z\vert}\le \frac{C}{\vert z\vert}
\end{multline*}
for $\vert z\vert < 1$. So, the lower estimate in \eqref{f_sol_bound} also holds.
\begin{flushright}
$\square$
\end{flushright}

\textbf{Proof of Lemma \ref{l_summand_order_1}:} By the definition of $r_a$ there exist $\nu_0\in\left\lbrace 1,\ldots,d_a\right\rbrace$ and $\tau_0\in\left\lbrace 1,\ldots,d_{r_a}\right\rbrace$ such that $\left(\mathfrak{S}_{a,r_a}\right)_{\nu_0,\tau_0}>0$.

\textbf{(1)} For $b<r_a$ and $\nu\in\left\lbrace 1,\ldots, d_a\right\rbrace$, $\tau\in\left\lbrace 1,\ldots, d_b\right\rbrace$ we have:\\
\begin{equation}
m^{(a)}_{\nu}(w_q)m^{(b)}_{\tau}(w_q) = o\left(m^{(a)}_{\nu_0}(w_q)m^{(r_a)}_{\tau_0}(w_q)\right) = o(1), \quad \text{as}\quad q\to\infty.
\label{f_oO}
\end{equation}
In \eqref{f_oO} the first equality follows from \eqref{f_diff_blocks} and the second from \eqref{f_F_entrywise}.

\textbf{(2)} By the definition of the block structure of the matrix $S$ we have that for any $\nu\in\left\lbrace 1,\ldots,d_a\right\rbrace$ and $\tau\in\left\lbrace 1,\ldots,d_{r_a}\right\rbrace$ it holds that\\
\begin{equation*}
m^{(a)}_{\nu}m^{(r_a)}_{\tau}\sim m^{(a)}_{\nu_0}m^{(r_a)}_{\tau_0} \quad \text{WRT} \quad \left\lbrace w_q\right\rbrace_{q=1}^\infty, 
\end{equation*}
So, it is sufficient to check that\\
\begin{equation}
m^{(a)}_{\nu_0}m^{(r_a)}_{\tau_0} \sim 1\quad \text{WRT} \quad \left\lbrace w_q\right\rbrace_{q=1}^\infty
\end{equation}
From (\ref{f_F_entrywise}) we have an upper estimate:\\
\begin{equation}
\left\vert m^{(a)}_{\nu_0}(w_q)\right\vert\cdot\left\vert m^{(r_a)}_{\tau_0}(w_q)\right\vert < \frac{1}{\left(\mathfrak{S}_{a,r_a}\right)_{\nu_0,\tau_0}}<C,\quad \forall l\in\mathbb{N}
\end{equation} 
Suppose that the lower estimate does not hold. Therefore, there exists a subsequence $\left\lbrace w_{q_u}\right\rbrace_{u=1}^\infty$ such that\\
\begin{equation}
\lim_{u\to\infty} \left\vert m^{(a)}_{\nu_0}(w_{q_u})\right\vert\cdot\left\vert m^{(r_a)}_{\tau_0}(w_{q_u})\right\vert = 0
\end{equation}
Hence, for all $\nu \in\left\lbrace 1,\ldots, d_a\right\rbrace ,\tau\in\left\lbrace 1,\ldots,d_{r_a}\right\rbrace$ it holds that\\
\begin{equation}
\lim_{u\to\infty} \left\vert m^{(a)}_\nu(w_{q_u})\right\vert\cdot\left\vert m^{(r_a)}_\tau(w_{q_u})\right\vert = 0
\label{f_summand_order1}
\end{equation} 
Note that for every $\nu\in\left\lbrace 1,\ldots,d_a\right\rbrace$ we have\\
\begin{equation}
\left\vert zm^{(a)}_\nu\right\vert = \left\vert m^{(r_a)}_1 m^{(a)}_\nu\right\vert \cdot\left\vert \frac{z}{m^{(r_a)}_1}\right\vert <C\left\vert m^{(a)}_\nu m^{(r_a)}_1\right\vert,
\end{equation}
here we used the lower estimate from Proposition \ref{prop_sol_bound}. Hence,\\
\begin{equation}
\lim_{u\to\infty} {w_{q_u}m^{(a)}_\nu(w_{q_u})}=0
\label{f_summand_order2}
\end{equation}
Therefore, from \eqref{f_summand_order1}, \eqref{f_summand_order2} and the first part of Lemma \ref{l_summand_order_1} it follows that for given $a$ all summands in the RHS of (\ref{nfVDE1}) tend to zero along the sequence $\left\lbrace w_{q_u}\right\rbrace_{u=1}^\infty$, while the LHS of (\ref{nfVDE1}) is a non-zero constant, which is a contradiction. So, the second part of Lemma \ref{l_summand_order_1} is also proved. 
\begin{flushright}
$\square$
\end{flushright}

\begin{lemma}\label{l_ineq_for_r}
Under Assumption \ref{assump_sol_ord} for all $a\in\left\lbrace 1,\ldots,p-1\right\rbrace$ it holds that\\
\begin{equation*}
r_{a+1}<r_a
\end{equation*}
\end{lemma}

\textbf{Proof of Lemma \ref{l_ineq_for_r}:} From Lemma \ref{l_summand_order_1} follows that\\
\begin{equation}
m^{(a)}_1m^{(r_a)}_1\sim 1 \quad \text{WRT} \quad \left\lbrace w_q\right\rbrace_{q=1}^\infty
\label{f_ineq_for_r1} 
\end{equation}
and\\
\begin{equation}
m^{(a+1)}_1m^{(r_{a+1})}_1\sim 1 \quad \text{WRT} \quad \left\lbrace w_q\right\rbrace_{q=1}^\infty 
\label{f_ineq_for_r2}
\end{equation}
Dividing (\ref{f_ineq_for_r1}) by (\ref{f_ineq_for_r2}) we get\\
\begin{equation}
\frac{m^{(a)}_1(w_q)}{m^{(a+1)}_1(w_q)}\sim \frac{m^{(r_{a+1})}_1(w_q)}{m^{(r_a)}_1(w_q)} \quad \text{WRT} \quad \left\lbrace w_q\right\rbrace_{q=1}^\infty
\label{f_ineq_for_r3}
\end{equation}
But the sequence in the LHS of the equivalence (\ref{f_ineq_for_r3}) tends to zero by \eqref{f_diff_blocks}, so\\
\begin{equation}
\lim_{q\to\infty}{\frac{m^{(r_{a+1})}_1(w_q)}{m^{(r_a)}_1(w_q)}} = 0 
\end{equation} 
Therefore, $r_{a+1}<r_a$ again by \eqref{f_diff_blocks}.\\

\begin{flushright}
$\square$
\end{flushright}

\begin{lemma} \label{l_eq_for_r}
Under Assumption \ref{assump_sol_ord} for all $a\in\left\lbrace 1,\ldots,p\right\rbrace$ it holds that\\
\begin{equation*}
r_a = p+1-a
\end{equation*}
\end{lemma}

\textbf{Proof of Lemma \ref{l_eq_for_r}:} The result immediately follows from the following two observations:\\ 
1) from Lemma \ref{l_ineq_for_r} we have $r_p<r_{p-1}<\ldots <r_1$ and\\
2) $1\le r_a \le p$ for all $a=1,\ldots ,p$.
\begin{flushright}
$\square$
\end{flushright}

\begin{lemma}\label{l_eq_d}
Under Assumptions \ref{assump_zm} and \ref{assump_sol_ord} for all $a\in\left\lbrace 1,\ldots,p\right\rbrace$ it holds that\\
\begin{equation*}
d_a=d_{p+1-a}
\end{equation*}
\end{lemma}

\textbf{Proof of Lemma \ref{l_eq_d}:} Fix $a_0\in\left\lbrace 1,\ldots,p\right\rbrace$. Summing up the equations of (\ref{nfVDE1}) for $a=a_0$ we obtain that\\
\begin{equation}
-d_{a_0} = \sum\limits_{\nu=1}^{d_{a_0}}zm^{(a_0)}_\nu + \sum_{b=1}^{p-a_0}\sum\limits_{\nu=1}^{d_{a_0}}\sum\limits_{\tau=1}^{d_b}\left(\mathfrak{S}_{a_0,b}\right)_{\nu,\tau} m^{(a_0)}_\nu m^{(b)}_\tau +\sum\limits_{\nu=1}^{d_{a_0}}\sum\limits_{\tau=1}^{d_{p+1-a_0}}\left(\mathfrak{S}_{a_0,p+1-a_0}\right)_{\nu,\tau} m^{(a_0)}_\nu m^{(p+1-a_0)}_\tau
\label{f_eq_d_1}
\end{equation}
Here by determining the limits of summation we used Lemma \ref{l_eq_for_r} for $a=a_0$.\\
Similarly for $a=p+1-a_0$:\\
\begin{multline}
-d_{p+1-a_0} = \sum\limits_{\nu=1}^{d_{p+1-a_0}}zm^{(p+1-a_0)}_\nu + \sum_{b=1}^{a_0-1}\sum\limits_{\nu=1}^{d_{p+1-a_0}}\sum\limits_{\tau=1}^{d_b}\left(\mathfrak{S}_{p+1-a_0,b}\right)_{\nu,\tau} m^{(p+1-a_0)}_\nu m^{(b)}_\tau +\\
+\sum\limits_{\nu =1}^{d_{p+1-a_0}}\sum\limits_{\tau =1}^{d_{a_0}}\left(\mathfrak{S}_{p+1-a_0,a_0}\right)_{\nu,\tau} m^{(p+1-a_0)}_\nu m^{(a_0)}_\tau
\label{f_eq_d_2}
\end{multline}
We subtract (\ref{f_eq_d_2}) from (\ref{f_eq_d_1}) and use that $\mathfrak{S}_{a_0,p+1-a_0} = \left(\mathfrak{S}_{p+1-a_0,a_0}\right)^T$: \\
\begin{multline}
d_{p+1-a_0}-d_{a_0} = \sum\limits_{\nu =1}^{d_{a_0}}zm^{(a_0)}_\nu - \sum\limits_{\nu=1}^{d_{p+1-a_0}}zm^{(p+1-a_0)}_\nu +\\
+ \sum_{b=1}^{p-a_0}\sum\limits_{\nu=1}^{d_{a_0}}\sum\limits_{\tau=1}^{d_b}\left(\mathfrak{S}_{a_0,b}\right)_{\nu,\tau} m^{(a_0)}_\nu m^{(b)}_\tau - \sum_{b=1}^{a_0-1}\sum\limits_{\nu =1}^{d_{p+1-a_0}}\sum\limits_{\tau =1}^{d_b}\left(\mathfrak{S}_{p+1-a_0,b}\right)_{\nu,\tau} m^{(p+1-a_0)}_\nu m^{(b)}_\tau
\label{f_eq_d_3}
\end{multline}
Now let us consider (\ref{f_eq_d_3}) for $z = w_q$, $q\in\mathbb{N}$. From the first statement of Lemma \ref{l_summand_order_1} it follows that every summand in the RHS of (\ref{f_eq_d_3}) tends to zero as $q\to \infty$. Therefore, $d_{p+1-a_0}=d_{a_0}$.\\
\begin{flushright}
$\square$
\end{flushright}

\textbf{Proof of Proposition \ref{prop_staircase}:} First two statements follow from Lemma \ref{l_eq_for_r} and the third from Lemma \ref{l_eq_d}.\\
\begin{flushright}
$\square$
\end{flushright}

\subsection{Behaviour of the solution of (\ref{VDE}) under Assumption \ref{assump_S}.} \label{subs_sol_ord}

While in the previous subsection Assumption \ref{assump_sol_ord} was imposed on the solution of VDE (\ref{VDE}), now we consider the matrix $S$ under Assumption \ref{assump_S} and prove that Assumption \ref{assump_sol_ord} holds. We still keep the technical Assumption \ref{assump_zm} that will be eliminated later in Section 5.\\

\begin{theorem}\label{t_sol_ord}
Under Assumptions \ref{assump_S} and \ref{assump_zm} there exists a subsequence $\left\lbrace z_{l_q}\right\rbrace_{q=1}^\infty$ of the sequence $\left\lbrace z_l\right\rbrace_{l=1}^\infty$ such that for every $k\in\left\lbrace 1,\ldots,n-1\right\rbrace$ either $m_k\sim m_{k+1}$ with respect to the sequence $\left\lbrace z_{l_q}\right\rbrace_{q=1}^{\infty}$ or\\
\begin{equation}
\lim_{q\to\infty} {\frac {m_k(z_{l_q})}{m_{k+1}(z_{l_q})}} = 0,
\label{f_sol_ord}
\end{equation} 
i.e. Assumption \ref{assump_sol_ord} holds.\\
\end{theorem}

From Theorem \ref{t_sol_ord} it follows that under Assumptions \ref{assump_S} and \ref{assump_zm} the matrix $S$ has a block staircase structure in the sense of Proposition \ref{prop_staircase}. Moreover, at the end of this section we will prove the following property of the anti-diagonal blocks of the matrix $S$:\\ 

\begin{prop}\label{prop_diag_irred}
Let the matrix $S$ satisfy Assumption \ref{assump_S}. Then for all $a\in\left\lbrace 1,\ldots,p\right\rbrace$ the matrix\\
\begin{equation*}
\left(\mathfrak{S}_{a,p+1-a}\right)\left(\mathfrak{S}_{a,p+1-a}\right)^T
\end{equation*}
is irreducible.\\
\end{prop}

The main idea of the proof of Theorem \ref{t_sol_ord} is the following lemma.

\begin{lemma}\label{l_comb}
Let $A$ be an $n\times n$ matrix with $0,1$ entries. Assume that $a_{i,j}=1$ for all $i,j\in\left\lbrace 1,\ldots,n\right\rbrace$ such that $i+j\in\left\lbrace n,n+1\right\rbrace$. Consider a permutation $\pi\in S_n$. Define an $n\times n$ matrix $B$ as follows:\\
\begin{equation}
b_{i,j}=a_{\pi(i),\pi(j)},\quad i,j\in\left\lbrace 1,\ldots, n\right\rbrace.
\end{equation}
Assume that for some fixed $k\in\left\lbrace 1,\ldots,n-1\right\rbrace$:\\
\begin{equation}
b_{i,j}=0 \quad \forall i\ge n-k+1, \quad \forall j\ge k+1.
\label{f_comb_assump_on_b}
\end{equation}
Then\\
\begin{equation}
\pi(i)\in\left\lbrace 1,\ldots,k\right\rbrace\quad \forall i\in\left \lbrace 1,\ldots,k\right\rbrace
\label{l_comb_1}
\end{equation}
and\\
\begin{equation}
\pi(i)\in\left\lbrace n+1-k,\ldots,n\right\rbrace\quad \forall i\in\left \lbrace n+1-k,\ldots,n\right\rbrace
\label{l_comb_2}
\end{equation}
\end{lemma}

\textbf{Proof of Lemma \ref{l_comb}:} It is sufficient to prove Lemma \ref{l_comb} in the case when\\
\begin{equation}
a_{i,j}=0 \quad \text{for}\quad i,j\in\left\lbrace 1,\ldots,n\right\rbrace\quad \text{such that} \quad i+j\notin \left\lbrace n,n+1\right\rbrace
\label{f_rem_comb}
\end{equation}
since the assumption in \eqref{f_comb_assump_on_b} becomes stronger if all $a_{i,j}$'s are set to zero for $i+j\notin\left\lbrace n,n+1\right\rbrace$. Thus we will assume that (\ref{f_rem_comb}) holds. Note, that after this assumption the matrices $A$ and $B$ are symmetric, so it is sufficient to prove the statement of Lemma \ref{l_comb} for\\
\begin{equation}
1\le k\le \left[\frac{n}{2}\right].
\label{f_rem_comb1}
\end{equation}
We also assume further that (\ref{f_rem_comb1}) holds.\\

We start with proving (\ref{l_comb_1}). The matrix $B$ can be represented in the following way\\
\begin{equation}
B = \begin{pmatrix}
B_{1,1} &B_{1,2} &B_{1,3}\\
B_{2,1} &B_{2,2} &0\\
B_{3,1} &0 &0,\\
\end{pmatrix}
\end{equation}
where $B_{1,1}$, $B_{1,2}$, $B_{1,3}$ and $B_{2,2}$ are $k\times k$, $k\times (n-2k)$, $k\times k$ and $(n-2k)\times (n-2k)$ matrices respectively, and\\
\begin{equation*}
B_{1,1}=B_{1,1}^T, \quad B_{2,2} = B_{2,2}^T,
\end{equation*}
\begin{equation}
B_{2,1}=B_{1,2}^T, \quad B_{3,1} = B_{1,3}^T.
\end{equation}
Consider the rows $\pi(n-k+1),\pi(n-k+2),\ldots,\pi(n)$ of the matrix $A$. There exists a set $K$ of $k$ columns such that each "1", which is contained in these rows, also belongs to one of the columns from the set $K$. But it is possible only for $k$ lowest rows of the matrix $A$, i.e.\\
\begin{equation}
\left\lbrace \pi(j)\right\rbrace_{j=n-k+1}^n \subset \left\lbrace n-k+1,\ldots,n\right\rbrace.
\end{equation}
Hence, (\ref{l_comb_1}) holds.\\

Each row of the matrix $A$ contains at least one "1". If for some $j\in\left\lbrace 1\ldots,k\right\rbrace$ it holds that $\pi(j)>k$, then from (\ref{l_comb_1}) it follows that $b_{\pi(j),\pi(n+1-j)}=1$. But $\pi(j)>k$ and $\pi(n-j+1)>n-k+1$ which is contradicting to \eqref{f_comb_assump_on_b}. So, (\ref{l_comb_2}) also holds.

\begin{flushright}
$\square$
\end{flushright}

\textbf{Proof of Theorem \ref{t_sol_ord}:} Extract a subsequence $\left\lbrace z_{l_q}\right\rbrace_{q=1}^{\infty}$ of the sequence $\left\lbrace z_l\right\rbrace_{l=1}^\infty$ such that for all $i,j\in\left\lbrace 1,\ldots,n\right\rbrace$ one of the following statements about the sequence $\left\lbrace\frac{\vert m_{i}(z_{l_q})\vert}{\vert m_{j}(z_{l_q})\vert}\right\rbrace_{q=1}^{\infty}$ holds:\\
\begin{enumerate}
\item{it is bounded;}\\
\item{it tends to zero;}\\
\item{it tends to infinity.}\\
\end{enumerate}
For any permutation $\pi\in S_n$ define the $n\times n$ matrix $S^{\pi}$ as follows: for $i,j\in\left\lbrace 1,\ldots,n\right\rbrace$\\
\begin{equation}
s^{\pi}_{i,j}=s_{\pi(i),\pi(j)}.
\end{equation}
Let $\pi_0\in S_n$ be a permutation such that for the solution of (\ref{VDE}) corresponding to the matrix $S^{\pi_0}$ Assumption \ref{assump_sol_ord} holds. Such a permutation exists but it is not necessarily unique. The block structure $\left(\mathfrak{S}^{\pi_0}_{a,b}\right)_{a,b}^n$ and the set $\left\lbrace d^{\pi_0}_a\right\rbrace_{a=1}^p$ are defined for the matrix $S^{\pi_0}$ as in Subsection \ref{subs_sol_ord}. From Lemma \ref{l_comb} for the permutation $\pi_0$ and $k:=d^{\pi_0}_1+\ldots +d^{\pi_0}_a$, $a\in\left\lbrace 1,\ldots,p-1\right\rbrace$ we obtain that\\ 
\begin{equation}
\pi_0\left(\left\lbrace 1,\ldots ,d^{\pi_0}_1+\ldots + d^{\pi_0}_a\right\rbrace\right) = \left\lbrace 1,\ldots ,d^{\pi_0}_1+\ldots + d^{\pi_0}_a\right\rbrace.  
\end{equation} 
Therefore,\\
\begin{equation*}
\pi_0\left(\left\lbrace 1,\ldots ,d^{\pi_0}_1\right\rbrace\right) = \left\lbrace 1,\ldots ,d^{\pi_0}_1\right\rbrace,
\end{equation*}
\begin{equation*}
\pi_0\left(\left\lbrace d^{\pi_0}_1+1,\ldots ,d^{\pi_0}_1+d^{\pi_0}_2\right\rbrace\right) = \left\lbrace d^{\pi_0}_1+1,\ldots ,d^{\pi_0}_1+d^{\pi_0}_2\right\rbrace,
\end{equation*}
\begin{equation*}
\ldots
\end{equation*}
\begin{equation*}
\pi_0\left(\left\lbrace d^{\pi_0}_1+\ldots + d^{\pi_0}_{p-1},\ldots ,d^{\pi_0}_1+\ldots + d^{\pi_0}_p\right\rbrace\right) = \left\lbrace d^{\pi_0}_1+\ldots+d^{\pi_0}_{p-1}+1,\ldots ,d^{\pi_0}_1+\ldots + d^{\pi_0}_p\right\rbrace.
\end{equation*}
It means that $\pi_0$ acts nontrivially only inside the blocks, but it leaves all blocks invariant, hence for every $k=1,\ldots,n$:\\
\begin{equation*}
m_{\pi_0(k)}\sim m_k \quad \text{WRT} \quad \left\lbrace z_{l_q}\right\rbrace_{q=1}^\infty.
\end{equation*}
Thus the solution of (\ref{VDE}) corresponding to the matrix $S$ also satisfies Assumption \ref{assump_sol_ord}. 
\begin{flushright}
$\square$
\end{flushright}

Since $\mathfrak{S}_{a,p+1-a}$ is an anti-diagonal block of $S$, directly from Theorem \ref{t_sol_ord} it follows:\\

\begin{lemma}\label{l_diag_block}
Let the matrix $S$ satisfy Assumption \ref{assump_S}. Then for $a\in\left\lbrace 1,\ldots,p\right\rbrace$ and $\nu,\tau\in\left\lbrace 1,\ldots,d_a\right\rbrace$ we have:\\
\begin{equation}
\left(\mathfrak{S}_{a,p+1-a}\right)_{\nu,\tau}>0\quad \text{if}\quad \nu+\tau\in\left\lbrace d_a,d_a+1\right\rbrace.
\end{equation}
\end{lemma}

Proposition \ref{prop_diag_irred} is an easy consequence of Lemma \ref{l_diag_block} and the fact that all matrix elements of the matrix $S$ are non-negative.

\section{Properties of phases of the solution of \eqref{VDE}.}
We work under Assumptions \ref{assump_S} and \ref{assump_zm}. We consider only $z=z_{l_q}$, $q\in\mathbb{N}$, where the sequence $\left\lbrace z_{l_q}\right\rbrace_{q=1}^\infty$ is from Theorem \ref{t_sol_ord} and use the same notation as in Section \ref{sec2}. Denote again $w_q:=z_{l_q}$, $q\in\mathbb{N}$. The main intermediate result of this section is as follows:\\

\begin{prop}\label{prop_eq_phases}
Under Assumptions \ref{assump_S} and \ref{assump_zm} for every $a\in\left\lbrace 1,\ldots,p\right\rbrace$ and $\nu,\tau\in\left\lbrace 1,\ldots,d_a\right\rbrace$ it holds that\\
\begin{equation}
\lim_{q\to\infty} \operatorname{arg}{\frac{m^{(a)}_\nu(w_q)}{m^{(a)}_\tau(w_q)}}=0.
\end{equation}
\end{prop}
Proposition \ref{prop_eq_phases} means that not only the moduli of the components of $m$ corresponding to the same block are comparable, but their phases tend to be the same.\\

The following lemma shows how the operator $F$ acts on the vector $\frac{m}{\vert m\vert}$.\\

\begin{lemma}\label{lF2} 
\begin{equation}
F^2\frac{m}{\vert m \vert} = \frac{m}{\vert m \vert} + \left(\overline{z}I-zF\right)\vert m\vert,
\label{f_F2}
\end{equation}
where $m=m(z)$ is the solution of (\ref{VDE}).
\end{lemma}

\begin{rem}\label{rem_F2}
Lemma \ref{lF2} does not use Assumptions \ref{assump_S} and \ref{assump_zm}.
\end{rem}

\textbf{Proof of Lemma \ref{lF2}:} By taking the imaginary and the real part of (\ref{VDE}) we obtain respectively:\\
\begin{equation}
\frac{\operatorname{Im}(m)}{\vert m\vert ^2} = \operatorname{Im}(z) + S\operatorname{Im}(m)
\Leftrightarrow F\frac{\operatorname{Im}(m)}{\vert m\vert} = \frac{\operatorname{Im}(m)}{\vert m\vert} - \vert m\vert\operatorname{Im} (z),
\label{Im}
\end{equation}
\begin{equation}
-\frac{\operatorname{Re}(m)}{\vert m\vert ^2} = \operatorname{Re}(z) + S\operatorname{Re}(m) 
\Leftrightarrow F\frac{\operatorname{Re}(m)}{\vert m\vert} = -\frac{\operatorname{Re}(m)}{\vert m\vert} - \vert m\vert\operatorname{Re}(z).
\label{Re}
\end{equation}
From (\ref{Im}) and (\ref{Re}) it follows:\\
\begin{multline}
F\frac{m}{\vert m\vert} = F\frac{\operatorname{Re}(m)}{\vert m\vert} + iF\frac{\operatorname{Im}(m)}{\vert m\vert} = \frac{-\operatorname{Re}(m)+i\operatorname{Im}(m)}{\vert m \vert} - \vert m\vert z = - \frac{\overline{m}}{\vert m\vert} - \vert m\vert z.
\label{Fm}
\end{multline}
Taking $F$ of both sides of (\ref{Fm}) we get:\\
\begin{multline}F^2\frac{m}{\vert m\vert} = -\overline{F\frac{m}{\vert m\vert}} - zF\vert m\vert = \frac{m}{\vert m\vert} +\vert m\vert \overline{z} -zF\vert z\vert = \frac{m}{\vert m\vert} + \left(\overline{z}I-zF\right)\vert m\vert. 
\label{F2m}
\end{multline}
\begin{flushright}
$\square$
\end{flushright}

\begin{cor} \label{cor_error}
Under Assumption \ref{assump_zm} it holds that\\
\begin{equation}
F^2(w_q)\frac{m(w_q)}{\vert m(w_q) \vert} = \frac{m(w_q)}{\vert m(w_q) \vert} + o(1),\quad q\rightarrow \infty
\label{f_error}
\end{equation}
\end{cor}

\textbf{Proof of Corollary \ref{cor_error}:} Let us find an upper bound for the last term of (\ref{f_F2}):\\
\begin{multline*} 
\|\left(\overline{w_q}I-w_qF(w_q)\right)\cdot\vert m(w_q)\vert\| \le \vert w_q\vert\cdot \|m(w_q)\| +\vert w_q\vert \cdot \|F(w_q)m(w_q)\| \le \\
\le \vert w_q\vert\cdot \left(1+\|F(w_q)\|\right)\cdot\|m(w_q)\|\le 2\vert w_q\vert\cdot\|m(w_q)\| \rightarrow 0 \quad \text{as}\quad q\to\infty 
\end{multline*}
By stating that $2\vert w_q\vert\cdot\|m(w_q)\| \rightarrow 0$ as $q\to\infty$ we used Assumption \ref{assump_zm}.\\ 
\begin{flushright}
$\square$
\end{flushright}

\begin{rem} We know that $m_j(z)\neq 0$ and $\operatorname{Im}{m_j(z)}>0$ for all $j=1,\ldots,n$, $z\in\mathbb{H}$.  Hence, an argument of a function $m_j$ is well defined, and we can choose its branch in such a way that\\ 
\begin{equation*}
\operatorname{arg}{m_j(z)}\in\left(0,\pi\right).
\end{equation*}
\end{rem}

Since all entries of the matrix $F(z)$ are bounded along the sequence $\left\lbrace w_q\right\rbrace_{q=1}^{\infty}$ and for all $j=1,\ldots,n$ holds $\operatorname{arg}{m_j(z)}\in (0,\pi)$,  we can extract a subsequence $\left\lbrace w_{q_u}\right\rbrace_{u=1}^\infty$ of the sequence $\left\lbrace w_q\right\rbrace_{q=1}^\infty$ such that\\
\begin{equation}
\exists \lim_{u\to\infty}\operatorname{arg}{m^{(a)}_\nu (w_{q_u})}=:\varphi^{(a)}_\nu,\quad \forall a\in \left\lbrace 1,\ldots,p\right\rbrace; \nu\in\left\lbrace 1,\ldots,d_a\right\rbrace
\end{equation}
and\\
\begin{equation}
\exists \lim_{u\to\infty}f_{j,k}\left(w_{q_u}\right) =:f^0_{j,k},\quad \forall j,k\in\left\lbrace 1,\ldots,n\right\rbrace,
\end{equation}
where $F(z) = \left(f_{j,k}(z)\right)_{j,k=1}^n$ and $F^0: = \left(f^0_{j,k}\right)_{j,k=1}^n$. For simplicity of the  notation we assume that the subsequence $\left\lbrace w_{q_u}\right\rbrace_{u=1}^\infty$ coincides with the sequence $\left\lbrace w_q\right\rbrace_{q=1}^\infty$.\\

Let us analyse the structure of the matrix $F^0$. The block structure $\left(\mathfrak{S}_{a,b}\right)_{a,b=1}^p$ of the  matrix $S$ induces the block structure in the matrix $F^0$: $\left(\mathfrak{F}^0_{a,b}\right)_{a,b=1}^p$, where $\mathfrak{F}^0_{a,b}$ and $\mathfrak{S}_{a,b}$ are $d_a\times d_b$ matrices.

\begin{lemma}\label{l_F0_cond}
The matrix $F^0$ meets the following conditions:\\
\begin{enumerate}
\item{$F^0$ is real-symmetric with non-negative entries.}\\
\item{$\|F^0\|\le 1$.}
\item{For $a,b\in\left\lbrace 1,\ldots,p\right\rbrace$ such that $a+b\neq p+1$ all entries of the matrix $\mathfrak{F}^0_{a,b}$ are zeros.}\\
\item{For $a\in\left\lbrace 1,\ldots,p\right\rbrace$ and $\nu,\tau\in\left\lbrace 1,\ldots,d_a\right\rbrace$:\\
\begin{equation*}
\mathfrak{F}^0_{\nu,\tau} \neq 0 \quad \text{if and only if} \quad \mathfrak{S}^0_{\nu,\tau} \neq 0.
\end{equation*}}\\
\end{enumerate}
\end{lemma}

Lemma \ref{l_F0_cond} immediately follows from the definition of the matrix $F^0$.\\

\textbf{Proof of Proposition \ref{prop_eq_phases}:} In order to prove that some sequence $\left\lbrace b_l\right\rbrace_{l=1}^\infty\subset \mathbb{C}$ converges to $b\in\mathbb{C}$ it is enough to show that for every subsequence $\left\lbrace b_{l_k}\right\rbrace_{k=1}^\infty$ of the sequence $\left\lbrace b_l\right\rbrace_{l=1}^\infty$ there exists a sub-subsequence $\left\lbrace b_{l_{k_j}}\right\rbrace_{j=1}^\infty$ converging to $b$. Using this fact for the sequence $\left\lbrace \operatorname{arg}{\frac{m^{(a)}_\nu(w_q)}{m^{(a)}_\tau(w_q)}}\right\rbrace_{q=1}^\infty$ we obtain that with the notation above it is sufficient to prove that for all $a\in\left\lbrace 1,\ldots,p\right\rbrace$ and $\nu,\tau\in\left\lbrace 1,\ldots,d_a\right\rbrace$ we have $\varphi^{(a)}_\nu = \varphi^{(a)}_\tau$.\\ 

Fix $a\in\left\lbrace 1,\ldots,p\right\rbrace$ and consider the canonical basis $\left\lbrace e_k\right\rbrace_{k=1}^n$ in $\mathbb{C}^n$. Denote\\
\begin{equation}
L_a = Lin\left\lbrace e_{d_1+\ldots + d_{a-1} + \nu}\right\rbrace_{\nu=1}^{d_a}
\end{equation}
Then $L_a$ is an invariant subspace for the operator $\left(F^0\right)^2$. Let $T$ be the restriction of $ \left(F^0\right)^2$ onto $L_a$. Projecting (\ref{f_error}) onto $L_a$ and letting $q$ go to infinity we obtain that $Tv_a = v_a$, where\\
\begin{equation}
v_a = \left(\exp{i\varphi^{(a)}_1},\ldots, \exp{i\varphi^{(a)}_{d_a}}\right).
\end{equation} 
Particularly, all components of the vector $v_a$ are non-zero. Let us note that in the canonical basis in $L_a$ the operator $T$ is given by the matrix\\
\begin{equation}
\left(\mathfrak{F}^0_{a,p+1-a}\right)\left(\mathfrak{F}^0_{a,p+1-a}\right)^T
\end{equation}
So, from Proposition \ref{prop_diag_irred} and Lemma \ref{l_F0_cond} it follows that the matrix $T$ is irreducible. Therefore, we already know about the matrix $T$ that it is irreducible with non-negative entries, $\|T\|\le 1$ and $Tv_a=v_a$, $v_a\neq 0$. Thus, $\|T\|=1$ and, applying Perron-Frobenius theorem [2, Theorem 8.4.4] to $T$, we obtain that $v_a=w_a\cdot\exp{i\varphi^{(a)}}$, where $\varphi^{(a)}\in\left[0,\pi\right]$ and $w_a$ is the Perron-Frobenius eigenvector of $T$. Hence for all $\nu,\tau\in\left\lbrace 1,\ldots,d_a\right\rbrace$ we have $\varphi^{(a)}_{\nu}=\varphi^{(a)}=\varphi^{(a)}_{\tau}$.
\begin{flushright}
$\square$
\end{flushright}

\section{Proof of the sequential version of Theorem \ref{t_main} for the absolute value.}
In this section we adhere Assumptions \ref{assump_S}, \ref{assump_zm} and prove Theorem \ref{t_main_seq}, which is a sequential version of Theorem \ref{t_main} for absolute value. We will replace the sequential approach by uniform one in Section 5 and in Section 6 Theorem \ref{t_main} will be proved in the full generality.\\

\begin{theorem}\label{t_main_seq}
Under Assumptions \ref{assump_S}, \ref{assump_zm} and with notations of Theorem \ref{t_sol_ord} there exist positive constants $\left\lbrace c_k\right\rbrace_{k=1}^n$ such that for each $k\in\left\lbrace 1,\ldots, n\right\rbrace$ we have\\
\begin{equation}
\lim_{q\to\infty} {\left(\vert m_k(w_q)\vert\cdot\vert w_q\vert ^{-\left(1-\frac{2k}{n+1}\right)}\right)} = c_k, 
\label{f_main_seq}
\end{equation}
where $w_q = z_{l_q}$ for $q\in\mathbb{N}$ and $\left\lbrace c_k\right\rbrace_{k=1}^n$ depend only on the entries of $S$, i.e. if Assumptions \ref{assump_zm} and \ref{assump_sol_ord} hold for another sequence $\left\lbrace w'_q\right\rbrace_{q=1}^\infty$, then the limits \eqref{f_main_seq} are the same.\\ 
\end{theorem}

First we prove the following Lemma:\\

\begin{lemma}\label{l_sys}
Under Assumptions \ref{assump_S} and \ref{assump_zm} we have:\\
\begin{enumerate}
\item{For $k\in\left\lbrace 1,\ldots,n\right\rbrace$: $m_km_{n+1-k}\sim 1$ WRT $\left\lbrace w_q\right\rbrace_{q=1}^\infty$.}\\
\item{$m_1m_{n-1}\sim zm_n$ WRT $\left\lbrace w_q\right\rbrace_{q=1}^\infty$.}\\
\item{For $k\in\left\lbrace 2,\ldots,n\right\rbrace$: $m_km_{n-k}\sim m_{n+1-k}m_{k-1}$ WRT $\left\lbrace w_q\right\rbrace_{q=1}^\infty$.}\\
\end{enumerate}
Particularly, denoting $m_0(z):=z$ and $m_{n+1}(z):=\frac{1}{z}$ for $k\in\left\lbrace 1,\ldots,n\right\rbrace$ we have\\
\begin{equation}
\frac{m_k}{m_{k+1}}\sim\frac{m_0}{m_1}\quad \text{WRT} \quad \left\lbrace w_q\right\rbrace_{q=1}^\infty.
\label{f_same_ratio}
\end{equation}
\end{lemma}

\textbf{Proof of Lemma \ref{l_sys}:} The first statement follows directly from Lemmas \ref{l_summand_order_1} and \ref{l_eq_d}. Now we start verifying the second one. By subtracting the last equation of (\ref{VDE}) from the first one we have\\
\begin{equation}
0 = \left(zm_1 + \sum\limits_{j=1}^{n-1} s_{1,j}m_1m_j\right) - zm_n.
\label{f_sys1}
\end{equation}
In order to prove the second statement of Lemma \ref{l_sys} it is sufficient to show that\\
\begin{equation}
\sum\limits_{j=1}^{n-1} s_{1,j}m_1m_j \sim m_1m_{n-1} \quad \text{WRT} \quad \left\lbrace w_q\right\rbrace_{q=1}^\infty.
\label{f_sys2}
\end{equation}
Indeed, then from (\ref{f_sys1}) and (\ref{f_sys2}) we will have that\\
\begin{equation}
zm_n = zm_1 + \sum\limits_{j=1}^{n-1} s_{1,j}m_1m_j \sim m_1m_{n-1} \quad \text{WRT} \quad \left\lbrace w_q\right\rbrace_{q=1}^\infty.
\end{equation}
So, let us split the sum $\sum\limits_{j=1}^{n-1} s_{1,j}m_1m_j$ into two:\\
\begin{equation}
\sum\limits_{j=1}^{n-1} s_{1,j}m_1m_j = \sum\limits_{j\in\Omega_1} s_{1,j}m_1m_j + \sum\limits_{j\in\Omega_2} s_{1,j}m_1m_j,
\end{equation}
where\\
\begin{equation*}
\Omega_1 = \left\lbrace j\in\left\lbrace 1,\ldots,n-1\right\rbrace\mid m_j\sim m_{n-1} \quad \text{WRT} \quad \left\lbrace w_q\right\rbrace_{q=1}^\infty\right\rbrace,
\end{equation*}
\begin{equation*}
\Omega_2 = \left\lbrace 1,\ldots,n-1\right\rbrace \backslash \Omega_1.
\end{equation*}
From the definition of the set $\Omega_2$ and (\ref{f_sol_ord}) follows that\\
\begin{equation}
\sum\limits_{j\in\Omega_2} s_{1,j}m_1(w_q)m_j(w_q) = o\left(m_1(w_q)m_{n-1}(w_q)\right),\quad\text{as}\quad q\to\infty. 
\label{f_sys3}
\end{equation}
Proposition \ref{prop_eq_phases} shows that the phases of summands of the sum $\sum\limits_{j\in\Omega_1} s_{1,j}m_1(w_q)m_j(w_q)$ become equal as $q$ tends to infinity. Therefore, there is no cancellation in the sum and\\
\begin{equation}
\sum\limits_{j\in\Omega_1} s_{1,j}m_1(w_q)m_j(w_q) \sim m_1(w_q)m_{n-1}(w_q),\quad\text{as}\quad q\to\infty.
\label{f_sys4}
\end{equation}
(\ref{f_sys2}) follows from (\ref{f_sys3}) and (\ref{f_sys4}), which completes the proof of the second statement.\\

The third statement of Lemma \ref{l_sys} can be proved similarly to the second one.\\

Now let us verify (\ref{f_same_ratio}). From the first three statements of Lemma \ref{l_sys} we obtain that\\
\begin{equation}
\frac{m_k}{m_{k+1}} \sim \frac{m_{n-k}}{m_{n-k+1}},\quad k\in\left\lbrace 0,1,\ldots,n\right\rbrace\quad \text{WRT} \quad \left\lbrace w_q\right\rbrace_{q=1}^\infty
\label{f_eq_ratio_1}
\end{equation}
and\\
\begin{equation}
\frac{m_{k-1}}{m_k} \sim \frac{m_{n-k}}{m_{n-k+1}},\quad k\in\left\lbrace 1,\ldots,n\right\rbrace\quad \text{WRT} \quad \left\lbrace w_q\right\rbrace_{q=1}^\infty.
\label{f_eq_ratio_2}
\end{equation}
So, for $k\in\left\lbrace 1,\ldots,n\right\rbrace$ using (\ref{f_eq_ratio_1}) and (\ref{f_eq_ratio_2}) we get\\
\begin{equation}
\frac{m_k}{m_{k+1}} \sim \frac{m_{n-k}}{m_{n-k+1}}\sim \frac{m_{k-1}}{m_k} \quad \text{WRT} \quad \left\lbrace w_q\right\rbrace_{q=1}^\infty,
\label{f_eq_ratio_3}
\end{equation}
which is equivalent to \eqref{f_same_ratio}.\\
\begin{flushright}
$\square$
\end{flushright}

\textbf{Proof of Theorem \ref{t_main_seq}:} We prove Theorem \ref{t_main_seq} in two steps: in the first part of the proof we show that\\
\begin{equation}
\vert m_k(w_q)\vert\cdot\vert w_q\vert ^{-\left(1-\frac{2k}{n+1}\right)}\sim 1 \quad \text{WRT}\quad \left\lbrace w_q\right\rbrace_{q=1}^\infty
\end{equation}
and in the second part we prove the existence of $\left\lbrace c_k\right\rbrace_{k=1}^n$.\\

\textbf{Part 1.} Multiplying all equivalences (\ref{f_same_ratio}) for $k=0,1,\ldots,n$ we get\\
\begin{equation*}
\frac{m_0}{m_1}\cdot\frac{m_1}{m_2}\cdots \frac{m_n}{m_{n+1}} \sim \left(\frac{m_0}{m_1}\right)^{n+1} \quad \text{WRT} \quad \left\lbrace w_q\right\rbrace_{q=1}^\infty.
\end{equation*}
So,\\
\begin{equation*}
\left(\frac{m_0}{m_1}\right)^{n+1}\sim z^2 \Rightarrow \frac{m_0}{m_1} \sim z^{\frac{2}{n+1}} \quad \text{WRT} \quad \left\lbrace w_q\right\rbrace_{q=1}^\infty.
\end{equation*}
Therefore, for every $k\in\left\lbrace 1,\ldots,n\right\rbrace$ we have:\\
\begin{equation}
m_k = \frac{m_k}{m_{k+1}}\cdot \frac{m_{k+1}}{m_{k+2}}\cdots \frac{m_n}{m_{n+1}}\cdot m_{n+1} \sim \left(\frac{m_0}{m_1}\right)^{n-k+1} \frac{1}{z}\sim \vert z\vert^{\frac{2(n-k+1)}{n+1}-1} = \vert z\vert^{1-\frac{2k}{n+1}} 
\end{equation}
with respect to the sequence $\left\lbrace w_q\right\rbrace_{q=1}^\infty$. It gives that there exist positive constants $c$ and $C$ such that for all $k\in\left\lbrace 1,\ldots,n\right\rbrace$ and $q\in\mathbb{N}$ we have\\
\begin{equation}
c<\vert m_k(w_q)\vert\cdot\vert w_q\vert ^{-\left(1-\frac{2k}{n+1}\right)}<C.
\end{equation}

\textbf{Part 2.} Consider an arbitrary subsequence $\left\lbrace w_{q_u}\right\rbrace_{u=1}^\infty$ of the sequence $\left\lbrace w_q\right\rbrace_{q=1}^\infty$ such that for $k\in\left\lbrace 1,\ldots,n\right\rbrace$ there exist the limits\\
\begin{equation}
\lim_{u\to\infty} {\left(\vert m_k(w_{q_u}\vert\cdot\vert w_{q_u}\vert ^{-\left(1-\frac{2k}{n+1}\right)}\right)} = \hat{c}_k.
\end{equation}
Similarly to the proof of Lemma \ref{l_sys} we obtain that $\left\lbrace \hat{c}_k\right\rbrace_{k=1}^n$ satisfy the following conditions:\\
\begin{equation}
\begin{cases}
s_{k,n+1-k}\hat{c}_k\hat{c}_{n+1-k}=1, &1\le k\le \left\lceil\frac{n}{2}\right\rceil\\
s_{1,n-1}\hat{c}_1\hat{c}_{n-1} = \hat{c}_n\\
s_{k,n-k}\hat{c}_k\hat{c}_{n-k} = s_{n+1-k,k-1}\hat{c}_{n+1-k}\hat{c}_{k-1}, &2\le k\le \left[\frac{n}{2}\right]\\
\end{cases}
\label{f_c_sys}
\end{equation}
Let us note that in \eqref{f_c_sys} the number of variables equals to the number of equations. After taking the logarithm of the equations of (\ref{f_c_sys}) we obtain that $\left\lbrace \log{\hat{c}_k}\right\rbrace_{k=1}^n$ solves the following linear system:\\
\begin{equation}
\begin{cases}
\log{\hat{c}_k}+\log{\hat{c}_{n+1-k}}=-\log{s_{k,n+1-k}}, &1\le k\le \left\lceil\frac{n}{2}\right\rceil\\
\log{\hat{c}_1}+\log{\hat{c}_{n-1}} - \log{\hat{c}_n} = -\log{s_{1,n-1}}\\
\log{\hat{c}_k}+ \log{\hat{c}_{n-k}} - \log{\hat{c}_{n+1-k}} - \log{\hat{c}_{k-1}}=\log{s_{n+1-k,k-1}}-\log{s_{k,n-k}}, &2\le k\le \left[\frac{n}{2}\right]\\
\end{cases}
\label{f_logc_sys}
\end{equation}
Note that $s_{j,k}>0$ for $j+k\in\left\lbrace n,n+1\right\rbrace$ by Assumption \ref{assump_S}. If we show that \eqref{f_logc_sys} uniquely defines $\left\lbrace \log{\hat{c}_k}\right\rbrace_{k=1}^n$, then we will automatically obtain that the limits $\left\lbrace \hat{c}_k\right\rbrace_{k=1}^n$ are the same for every convergent subsequence of the sequence $\left\lbrace w_q\right\rbrace_{q=1}^\infty$ and depend only on the entries of $S$, which is equivalent to the statement of Theorem \ref{t_main_seq}.\\

So, in order to finish the proof of Theorem \ref{t_main_seq} it is sufficient to show that the homogeneous linear system corresponding to \eqref{f_logc_sys} has only zero solution. Consider a system\\
\begin{equation}
\begin{cases}
x_k + x_{n+1-k}=0, &1\le k\le \left\lceil\frac{n}{2}\right\rceil\\
x_1 + x_{n-1} = x_n\\
x_k + x_{n-k} = x_{n+1-k} + x_{k-1}, &2\le k\le \left[\frac{n}{2}\right]\\
\end{cases}
\label{f_x_sys}
\end{equation}
We need to show that $x_1=x_2=\ldots=x_n=0$. Define $x_0:=0$ and $x_{n+1}:=0$. From \eqref{f_x_sys} it follows that\\
\begin{equation}
\begin{cases}
x_k + x_{n+1-k} = 0, &0\le k\le n+1\\
x_k + x_{n-k} = x_{n+1-k} + x_{k-1}, &1\le k\le n\\
\end{cases}
\end{equation}
After the simple arithmetic transformations similar to \eqref{f_eq_ratio_1}, \eqref{f_eq_ratio_2} and \eqref{f_eq_ratio_3} we obtain that\\
\begin{equation}
x_0-x_1=x_1-x_2=\ldots=x_n-x_{n+1}=\frac{x_0-x_{n+1}}{n+1}=0.
\end{equation}
Therefore, $x_k=0$ for all $k\in\left\lbrace 1,\ldots,n\right\rbrace$.\\
\begin{flushright}
$\square$
\end{flushright}

\section{Elimination of the technical Assumption \ref{assump_zm} and verification of Theorem \ref{t_main} for the absolute value.}

In this section we show that Assumption \ref{assump_S} implies Assumption \ref{assump_zm} and prove Theorem \ref{t_main_abs_val} below which is the absolute value version of Theorem \ref{t_main}.\\

\begin{prop}\label{prop_zm}
Let the matrix $S$ satisfy Assumption \ref{assump_S}. Then for all $k\in\left\lbrace 1,\ldots,n\right\rbrace$ $zm_k(z)$ tends to zero as $z$ tends to zero from complex upper half-plane, i.e. Assumption \ref{assump_zm} holds.\\
\end{prop}

\begin{theorem}\label{t_main_abs_val}
Under Assumption \ref{assump_S} there exist positive constants $\left\lbrace c_k\right\rbrace_{k=1}^n$ such that for all $k\in\left\lbrace 1,\ldots,n\right\rbrace$ we have\\
\begin{equation}
\lim_{z\to 0} {\left(\vert m_k(z)\vert \cdot\vert z\vert^{-\left(1-\frac{2k}{n+1}\right)}\right)} = c_k.
\label{f_main_abs_val}
\end{equation}
\end{theorem}

For a function $f:\mathbb{H}\rightarrow \mathbb{R}$ define:\\
\begin{equation}
\Lambda\left\lbrace f\right\rbrace := \left\lbrace w\in\overline{\mathbb{R}}\mid \exists \left\lbrace z_l\right\rbrace_{l=1}^{\infty}\subset\mathbb{H}, \lim_{l\to\infty} z_l=0: w=\lim_{l\to\infty} f\left(z_l\right)\right\rbrace,
\label{def_set_part_lim}
\end{equation}
i.e., $\Lambda\left\lbrace f\right\rbrace$ is a set of partial limits of a function $f$ as $z$ tends to zero from the complex upper half-plane. In (\ref{def_set_part_lim}) by $\overline{\mathbb{R}}$ we mean the following:\\
\begin{equation}
\overline{\mathbb{R}}:= \mathbb{R} \sqcup \left\lbrace \pm\infty\right\rbrace.
\end{equation}

\begin{lemma}\label{l_limit_set}
Let $f:\mathbb{H}\rightarrow\mathbb{R}$ be a continuous function. Then $\Lambda\left\lbrace f\right\rbrace$ is a connected set.
\end{lemma}

\textbf{Proof of Lemma \ref{l_limit_set}:} Let $a$ and $b$ be arbitrary elements of $\Lambda\left\lbrace f\right\rbrace$, $a<b$. Consider arbitrary $c\in(a,b)$. It is sufficient to show that $c\in\Lambda\left\lbrace f\right\rbrace$. As $a,b\in\Lambda\left\lbrace f\right\rbrace$, there exist sequences $\left\lbrace a_k\right\rbrace_{k=1}^{\infty}\subset \mathbb{H}$ and $\left\lbrace b_k\right\rbrace_{k=1}^{\infty}\subset \mathbb{H}$ such that\\
\begin{equation}
\lim_{k\to\infty}{f(a_k)} = a,\quad \lim_{k\to\infty}{f(b_k)} = b.
\end{equation}
We can find $K\in\mathbb{N}$ such that for every $k>K$ $f(a_k)<c<f(b_k)$. From the continuity of the function $f$ follows that for each $k>K$ there exists $c_k$ in the segment with endpoints in points $a_k$ and $b_k$ such that $f(c_k)=c$. It is clear that $c_k$ tends to zero as $k$ tends to infinity. Therefore, $c\in\Lambda\left\lbrace f\right\rbrace$.
\begin{flushright}
$\square$
\end{flushright}

\begin{lemma}\label{l_on_axes_ineq}
Let the matrix $S$ satisfy Assumption \ref{assump_S}. Then for all $j\in\left\lbrace 1,\ldots, n\right\rbrace$ $iy\cdot m_j(iy)$ tends to zero as $y$ tends to $+0$.\\ 
\end{lemma}

\textbf{Proof of Lemma \ref{l_on_axes_ineq}:} Assume the opposite, i.e. that $iy\cdot\|m(iy)\|$ does not tend to zero as $y$ tends to $+0$. Then there exists a sequence $\left\lbrace y_l\right\rbrace_{l=1}^\infty\subset\left\lbrace y\in\mathbb{R}\mid y>0\right\rbrace$ converging to zero such that\\
\begin{equation}
y_l\|m(iy_l)\|>c, \quad \forall l\in\mathbb{N}.
\label{f_not_converge}
\end{equation}
Up to extracting a subsequence we can assume that for all $j,k\in\left\lbrace 1,\ldots,n\right\rbrace$ one of the following statements about the sequence $\left\lbrace \frac{m_j(iy_l)}{m_k(iy_l)}\right\rbrace_{l=1}^\infty$ holds:\\
\begin{enumerate}
\item{$c<\left\vert\frac{m_j(iy_l)}{m_k(iy_l)}\right\vert<C$ for all $l\in\mathbb{N}$.}\\
\item{It tends to zero.}\\
\item{It tends to infinity.}\\
\end{enumerate}

There exists the permutation $\pi\in S_n$ such that for the matrix $S^\pi$ Assumption \ref{assump_sol_ord} holds. We define the block structure of $S^\pi$ in the same way as it was done in Subsection \ref{subs_staircase}. Note that $S^\pi$ has the block staircase structure in the sense of Proposition \ref{prop_staircase}. Let the function $m^\pi:\mathbb{H}\rightarrow\mathbb{C}^n$ solve \eqref{VDE} corresponding to the matrix $S^\pi$. Particularly, we have $m^\pi=\left(m_{\pi(1)},m_{\pi(2)},\ldots,m_{\pi(n)}\right)$. As in Subsection \ref{subs_staircase} we denote the corresponding partition of the vector $m^\pi$ into blocks by $m^\pi=\left(m^{(1),\pi},\ldots,m^{(p),\pi}\right)$, where $m^{(k),\pi}\in\mathbb{C}^{d_k}$. From \eqref{f_not_converge} and the definition of the block structure of $m^\pi$ it follows that\\ 
\begin{equation}
\|m^{(p),\pi}(z)\|\sim \frac{1}{z},\quad \text{WRT}\quad \left\lbrace iy_l\right\rbrace_{l=1}^\infty.
\label{f_as_1_over_z}
\end{equation}
Lemma \ref{l_summand_order_1} and \eqref{f_as_1_over_z} give that\\ 
\begin{equation}
\|m^{(1),\pi}(z)\|\sim z\quad \text{WRT}\quad \left\lbrace iy_l\right\rbrace_{l=1}^\infty.
\end{equation}

We will consider two cases: $d_1<d_p$ and $d_1\ge d_p$.\\

\textbf{(1)} Assume that $d_1<d_p$. From Proposition \ref{prop_staircase} it follows that for all $d_1<j\le n$ and $n-d_p+1\le k\le n$ it holds that $s_{j,k}=0$. Therefore, rows with indices $d_1<j\le n$ and columns with indices $n-d_p+1\le k\le n$ form a zero block with a perimeter $2(d_p+(n-d_1))\ge 2(n+1)$, which contradicts to the fact that all anti-diagonal entries of the matrix $S$ are positive.\\

\textbf{(2)} Assume that $d_1\ge d_p$. Note that the solution of \eqref{VDE} is pure imaginary when $z$ is pure imaginary. So, now \eqref{VDE} takes the following form:\\
\begin{equation}
1 = y_l\left\vert m^{(a),\pi}_\nu(y_l)\right\vert + \sum\limits_{b=1}^{n+1-a} \sum\limits_{\tau=1}^{d_j}\left(\mathfrak{S}^\pi_{a,b}\right)_{\nu,\tau} \left\vert m^{(a),\pi}_\nu(y_l)\right\vert\cdot\left\vert m^{(b),\pi}_\tau(y_l)\right\vert,
\label{f_case_2}
\end{equation}
for all $  a\in\left\lbrace 1,\ldots,p\right\rbrace$, $\nu\in\left\lbrace 1,\ldots,d_a\right\rbrace$ and $l\in\mathbb{N}$. Subtracting the sum of all equations of \eqref{f_case_2} corresponding to $a=p$ from  the sum of all equations of \eqref{f_case_2} corresponding to $a=1$ we obtain that\\
\begin{equation}
d_1-d_p = -\sum_{\nu=1}^{d_p}y_l\left\vert m^{(p),\pi}_\nu(y_l)\right\vert + \left(y_l\left\vert m^{(1),\pi}_\nu(y_l)\right\vert + \sum\limits_{b=1}^{n-1} \sum\limits_{\tau=1}^{d_j}\left(\mathfrak{S}^\pi_{1,b}\right)_{\nu,\tau} \left\vert m^{(1),\pi}_\nu(y_l)\right\vert\cdot\left\vert m^{(b),\pi}_\tau(y_l)\right\vert\right).
\label{f_case_21}
\end{equation}
From Lemma \ref{l_summand_order_1} it follows that the second sum in the RHS of \eqref{f_case_21} tends to zero as $l$ tends to infinity. Also from \eqref{f_as_1_over_z} we know that $ \sum_{\nu=1}^{d_p}y_l\left\vert m^{(p),\pi}_\nu(y_l)\right\vert>c$. Therefore, from \eqref{f_case_21} we have\\
\begin{equation}
d_1-d_p<-c+o(1), \quad l\to\infty
\end{equation}
contradicting to $d_1\ge d_p$, which finishes the proof of Lemma \ref{l_on_axes_ineq}.\\
\begin{flushright}
$\square$
\end{flushright}

\begin{cor}\label{cor_on_im_axes}
Let the matrix $S$ satisfy Assumption \ref{assump_S}. Then\\
\begin{equation}
\Lambda\left\lbrace \frac{\log{\|m\|}}{\log{\vert z\vert}}\right\rbrace \cap (-1,1) = \left\lbrace -\frac{n-1}{n+1}\right\rbrace.
\end{equation}
\end{cor}

\textbf{Proof of Corollary \ref{cor_on_im_axes}:} Consider a sequence of positive numbers $\left\lbrace y_l\right\rbrace_{l=1}^\infty$ converging to zero. Lemma \ref{l_on_axes_ineq} shows that Assumption \ref{assump_zm} holds for the equation (\ref{VDE}) and the sequence $\left\lbrace iy_l\right\rbrace_{l=1}^\infty$. Thus, applying Theorem \ref{t_main_seq} to this sequence we get that for some subsequence $\left\lbrace iy_{l_q}\right\rbrace_{q=1}^\infty$ of the sequence $\left\lbrace y_l\right\rbrace_{l=1}^\infty$ it holds that\\ 
\begin{equation}
\|m\| \sim \vert z\vert ^{-\frac{n-1}{n+1}}\quad \text{WRT} \quad \left\lbrace iy_{l_q}\right\rbrace_{q=1}^\infty.
\end{equation}
So,\\
\begin{equation}
-\frac{n-1}{n+1} \in \Lambda\left\lbrace \frac{\log{\|m\|}}{\log{\vert z\vert}}\right\rbrace \cap (-1,1).
\end{equation}
On the other hand, if $a\in \Lambda\left\lbrace \frac{\log{\|m\|}}{\log{\vert z\vert}}\right\rbrace \cap (-1,1)$, then there exists a sequence $\left\lbrace z_l\right\rbrace_{l=1}^\infty\subset \mathbb{H}$ converging to zero such that\\
\begin{equation}
\lim_{l\to\infty} \frac{\log{\|m\|}}{\log{\vert z\vert}} = a>-1.
\end{equation}
Particularly, Assumption \ref{assump_zm} holds for (\ref{VDE}) and the sequence $\left\lbrace z_l\right\rbrace_{l=1}^\infty$. Hence from Theorem \ref{t_main_seq} applied to the sequence $\left\lbrace z_l\right\rbrace_{l=1}^\infty$ it follows that $a=-\frac{n-1}{n+1}$.\\
\begin{flushright}
$\square$
\end{flushright}

\textbf{Proof of Proposition \ref{prop_zm}:} Assume the opposite, i.e. $z\|m(z)\|$ does not tend to zero as $z$ tends to zero from upper half-plane. It means that there exists a point $b\in \Lambda\left\lbrace \frac{\log{\|m\|}}{\log{\vert z\vert}}\right\rbrace$ ($b\in\overline{\mathbb{R}}$) such that $\vert b\vert\ge 1$. From Lemma \ref{l_limit_set} we know that the set $\Lambda\left\lbrace \frac{\log{\|m\|}}{\log{\vert z\vert}}\right\rbrace$ is connected and $-\frac{n-1}{n+1} \in \Lambda\left\lbrace \frac{\log{\|m\|}}{\log{\vert z\vert}}\right\rbrace$ from Corollary \ref{cor_on_im_axes}. Therefore, there are infinitely many points in the set $\Lambda\left\lbrace \frac{\log{\|m\|}}{\log{\vert z\vert}}\right\rbrace \cap (-1,1)$, which contradicts to Corollary \ref{cor_on_im_axes}.

\begin{flushright}
$\square$
\end{flushright}

\textbf{Proof of Theorem \ref{t_main}:} From Proposition \ref{prop_zm} and Theorem \ref{t_main_seq} it follows that for each sequence $\left\lbrace z_l\right\rbrace_{l=1}^\infty\subset\mathbb{H}$ converging to zero there exists a subsequence $\left\lbrace z_{l_q}\right\rbrace_{q=1}^\infty$ and positive constants $\left\lbrace c_k\right\rbrace_{k=1}^n$ such that for all $k\in\left\lbrace 1,\ldots,n\right\rbrace$ we have\\
\begin{equation*}
\lim_{q\to\infty} {\left(\vert m_k(z_{l_q})\vert \cdot\vert z_{l_q}\vert^{-\left(1-\frac{2k}{n+1}\right)}\right)} = c_k,
\end{equation*}
where $\left\lbrace c_k\right\rbrace_{k=1}^n$ depend only on $S$. Therefore, (\ref{f_main_abs_val}) holds.\\
\begin{flushright}
$\square$
\end{flushright}

\section{Verification of \eqref{f_main} for the argument and completion of the proof of Theorem \ref{t_main}.}

In this section we work under Assumption \ref{assump_S} and prove Theorem \ref{t_main_ph} below which is the argument version of Theorem \ref{t_main}. We also finish the proof of Theorem \ref{t_main}.\\

\begin{theorem}\label{t_main_ph}
Under Assumption \ref{assump_S} for all $k\in\left\lbrace 1,\ldots,n\right\rbrace$ we have\\
\begin{equation}
\lim_{z\to 0} {\left( \operatorname{arg}{m_k(z)} - \left(1-\frac{2k}{n+1}\right)\operatorname{arg}{z} \right)} = \frac{\pi k}{n+1},
\label{f_main_ph}
\end{equation}
where $z$ tends to zero from the complex upper half-plane.\\ 
\end{theorem}

The proof below is similar to the proof of Theorem \ref{t_main_seq} except some technical details. 

Define the set of partial limits, as $z$ tends to zero, of a vector-function $f=\left(f_1,\ldots,f_n\right):\mathbb{H}\rightarrow \mathbb{R}^n$ as follows:\\
\begin{multline}
\Lambda\left\lbrace \left(f_1,\ldots,f_n\right)\right\rbrace := \left\lbrace \left(w_1,\ldots,w_n\right)\in\overline{\mathbb{R}}^n\mid \exists \left\lbrace z_l\right\rbrace_{l=1}^{\infty}\subset\mathbb{H}, \lim_{l\to\infty} z_l=0,\forall j=\overline{1,n}: w_j=\lim_{l\to\infty} f_j\left(z_l\right) \right\rbrace .
\label{def_vect_part_lim}
\end{multline}
In this notation \eqref{f_main_ph} takes the following form\\
\begin{equation}
\Lambda_0:= \Lambda \left\lbrace \left(\operatorname{arg}{m_k(z)} - \left(1-\frac{2k}{n+1}\right)\operatorname{arg}{z}\right)_{k=1}^n\right\rbrace = \left\lbrace \left(\frac{2\pi k}{n+1}\right)_{k=1}^n\right\rbrace.
\end{equation}

\begin{lemma}\label{l_part_lim_vec}
Let $f$ be a continuous vector-function, $f:\mathbb{H}\rightarrow \mathbb{R}^n$, which is bounded in some neighbourhood of zero in $\mathbb{H}$. Then either $\Lambda\left\lbrace f\right\rbrace$ consists of a single point or $\Lambda\left\lbrace f\right\rbrace$ is infinite.
\end{lemma}

\textbf{Proof of Lemma \ref{l_part_lim_vec}:} Let us assume the opposite: $f$ has at least two partial limits in zero, but the set $A=\Lambda\left\lbrace f\right\rbrace$ is finite. Consider an arbitrary point $a\in A$. As $A$ is finite, there exists $r>0$ such that $D_r(a)\cap A = \left\lbrace a\right\rbrace$, where\\
\begin{equation*}
D_r(a) = \left\lbrace x\in\mathbb{R}^n \mid \|x-a\|_2\le r\right\rbrace
\end{equation*} 
is a closed ball in $\mathbb{R}^n$. Let $b$ be an element of $\Lambda\left\lbrace f\right\rbrace$ different from $a$. From the definition \eqref{def_vect_part_lim} it follows that there exist sequences $\left\lbrace a_k\right\rbrace_{k=1}^{\infty}\subset\mathbb{H}$ and $ \left\lbrace b_k\right\rbrace_{k=1}^{\infty}\subset\mathbb{H}$ converging to zero such that\\
\begin{equation*}
\lim_{k\to \infty}f(a_k) = a,\quad \lim_{k\to \infty}f(b_k) = b.
\end{equation*}
For some $K\in \mathbb{N}$ it holds that $f(a_k)\in D_r(a)$, $f(b_k)\notin D_k(a)$ for every $k\ge K$. Let us consider a segment with endpoints in $a_k$ and $b_k$, and note that $f$ maps it continuously into some continuous curve in $\mathbb{R}^n$ with endpoints $f(a_k)$ and $f(b_k)$. So, for every $k\ge K$ there exists $\lambda_k\in(0,1)$ such that $f(\lambda_ka_k+\left(1-\lambda_k\right)b_k)\in S_r(a)$, where\\
\begin{equation*}
S_r(a) = \left\lbrace x\in\mathbb{R}^n \mid \|x-a\|_2=r\right\rbrace
\end{equation*}
is a sphere in $\mathbb{R}^n$ of a radius $r$ with a center $a$. Denote $c_k:=\lambda_ka_k+\left(1-\lambda_k\right)b_k$. Obviously,\\
\begin{equation*}
\lim_{k\to\infty}c_k=0\quad \text{and} \quad f\left(\left\lbrace c_k\right\rbrace_{k=1}^{\infty}\right)\subset S_r(a).
\end{equation*}
Since $S_r(a)$ is a compact set, we can extract a subsequence $\left\lbrace c_{k_l}\right\rbrace_{l=1}^{\infty}$ of $\left\lbrace c_k\right\rbrace_{k=1}^{\infty}$ such that $f\left(c_{k_l}\right)$ converges to some point of $S_r(a)$. Hence,\\
\begin{equation*}
\Lambda\left\lbrace f\right\rbrace \cap S_r(a)\neq \emptyset,
\end{equation*}
which contradicts to the assumption that $D_r(a)\cap A = \left\lbrace a\right\rbrace$.
\begin{flushright}
$\square$
\end{flushright}

Denote $m_0(z):=z$. The following lemma describes the properties of the points of the set $\Lambda\left\lbrace \left(\operatorname{arg}{m_k(z)}\right)_{k=0}^n\right\rbrace$:\\

\begin{lemma}\label{l_sys_ph}
For every element $\left(\varphi_k\right)_{k=0}^n \in \Lambda\left\lbrace \left(\operatorname{arg}{m_k(z)}\right)_{k=0}^n\right\rbrace$ under Assumption \ref{assump_S} we have:\\
\begin{enumerate}
\item{For $k\in\left\lbrace 1,\ldots,n\right\rbrace$: $e^{i\varphi_k}e^{i\varphi_{n+1-k}}= -1$.}\\
\item{For $k\in\left\lbrace 1,\ldots,n\right\rbrace$: $e^{i\varphi_k}e^{i\varphi_{n-k}}=e^{i\varphi_{n+1-k}}e^{i\varphi_{k-1}}$.}\\
\end{enumerate}
Particularly, denoting $\varphi_{n+1}:=\pi-\varphi_0$ for $k\in\left\lbrace1,\ldots,n\right\rbrace$ we have\\
\begin{equation}
\frac{e^{i\varphi_k}}{e^{i\varphi_{k+1}}}=\frac{e^{i\varphi_0}}{e^{i\varphi_1}}.
\label{f_same_ratio_ph}
\end{equation}
\end{lemma}

The proof of Lemma \ref{l_sys_ph} is completely analogous to the proof of Lemma \ref{l_sys} and hence is omitted.\\

\textbf{Proof of Theorem \ref{t_main_ph}:} At first we will show that $\left(\frac{2\pi k}{n+1}\right)_{k=1}^n\in\Lambda_0$ and then verify the finiteness of the set $\Lambda_0$. 

\textbf{Part 1.} Consider $\left(\varphi_k\right)_{k=0}^n \in \Lambda\left\lbrace \left(\operatorname{arg}{m_k(z)}\right)_{k=0}^n\right\rbrace$, where $\varphi_0\in(0,\pi)$, e.g. $z$ tends to zero along the ray $\operatorname{arg}{z}=\varphi_0$. We want to prove that\\
\begin{equation}
\varphi_k =\frac{\pi k}{n+1}+ \left(1-\frac{2k}{n+1}\right)\varphi_0,
\label{f_expl_ph}
\end{equation}
which is equivalent to\\
\begin{equation}
\varphi_{n+1}-\varphi_n=\varphi_n-\varphi_{n-1}=\ldots=\varphi_1-\varphi_0.
\label{f_eq_dist_phi}
\end{equation}
From \eqref{f_same_ratio_ph} it follows that for every $k\in\left\lbrace 1,\ldots,n\right\rbrace$ there exists an integer $t_k$ such that\\
\begin{equation}
\varphi_{k+1}-\varphi_k = \varphi_1-\varphi_0+2\pi t_k,
\end{equation}
i.e. $\varphi_{k+1}+\varphi_0-\varphi_k-\varphi_1 = 2\pi t_k$. As $\varphi_j\in [0,\pi]$ for $j\in\left\lbrace 0,\ldots,n+1\right\rbrace$, we have that\\
\begin{equation}
\varphi_{k+1}+\varphi_0-\varphi_k-\varphi_1 \in [-2\pi,2\pi].
\end{equation} 
Particularly, for every $k\in\left\lbrace 1,\ldots,n\right\rbrace$ it holds that $t_k\in\left\lbrace -1,0,1\right\rbrace$  and\\
\begin{equation*}
t_k=-1\quad \Rightarrow\quad \varphi_0=\varphi_{k+1}=0,\quad \varphi_1=\varphi_k=\pi, 
\end{equation*}
\begin{equation*}
t_k=1\quad \Rightarrow\quad \varphi_0=\varphi_{k+1}=\pi,\quad \varphi_1=\varphi_k=0. 
\end{equation*}
But $\varphi_0\in(0,\pi)$, so all $t_k=0$ and \eqref{f_eq_dist_phi} holds. Therefore, it also holds that 
\begin{equation}
\left(\frac{2\pi k}{n+1}\right)_{k=1}^n\in\Lambda_0.
\label{f_lambda_incl}
\end{equation}

\textbf{Part 2.}  Fix a point $\left(\psi_1,\ldots,\psi_n\right)\in\Lambda_0$. From \eqref{def_vect_part_lim} we have that there exists a sequence $\left\lbrace z_l\right\rbrace_{l=1}^\infty\subset \mathbb{H}$ converging to zero such that\\
\begin{equation}
\psi_k := \lim_{l\to\infty} {\left( \operatorname{arg}{m_k(z_l)} - \left(1-\frac{2k}{n+1}\right)\operatorname{arg}{z_l} \right)},\quad k\in\left\lbrace 1,\ldots,n\right\rbrace,
\label{def_psi}
\end{equation}
in particular, the limits in \eqref{def_psi} exist. Up to extracting a subsequence we can assume that there exist limits\\
\begin{equation}
\varphi_k:=\lim_{l\to\infty}{\operatorname{arg}{m_k(z_l)}}, \quad k\in\left\lbrace 1,\ldots,n\right\rbrace.
\label{def_phi}
\end{equation}
As above we denote $\varphi_{n+1}:=\pi-\varphi_0$. So, from the definitions \eqref{def_psi} and \eqref{def_phi} we obtain that\\
\begin{equation}
\psi_k = \varphi_k - \left(1-\frac{2k}{n+1}\right)\varphi_0.
\label{f_psi_phi}
\end{equation}
Recall that $\left\lbrace \varphi_k\right\rbrace_{k=1}^{n+1}$ solves the system\\
\begin{equation}
\begin{cases}
\varphi_{k+1}-\varphi_k = \varphi_1-\varphi_0+2\pi t_k, &k\in\left\lbrace 1,\ldots,n\right\rbrace,\\
\varphi_{n+1}=\pi-\varphi_0,\\
\end{cases}
\label{sys_phi}
\end{equation}
where $t_k\in\left\lbrace -1,0,1\right\rbrace$ for $k\in\left\lbrace 1,\ldots,n\right\rbrace$. We will find $\left\lbrace \varphi_k\right\rbrace_{k=1}^{n+1}$ explicitly from \eqref{sys_phi}:\\
\begin{multline*}
\pi-2\varphi_0=\varphi_{n+1}-\varphi_0=(\varphi_{n+1}-\varphi_n)+(\varphi_n-\varphi_{n-1})+\ldots+(\varphi_1-\varphi_0)=\\
=(n+1)(\varphi_1-\varphi_0) + 2\pi\sum\limits_{k=1}^nt_k.
\end{multline*}
So for $k\in\left\lbrace 1,\ldots,n\right\rbrace$ we have\\
\begin{multline*}
\varphi_k = (\varphi_k-\varphi_{k-1})+(\varphi_{k-1}-\varphi_{k-2}+\ldots+(\varphi_1-\varphi_0)+\varphi_0=\varphi_0+k(\varphi_1-\varphi_0)+\sum\limits_{j=1}^kt_j =\\
= \left(1 - \frac{2k}{n+1}\right)\varphi_0 +\left(\frac{\pi k}{n+1} - \frac{2\pi k}{n+1}\sum\limits_{j=1}^n t_j + \sum\limits_{j=1}^k t_k\right). 
\end{multline*}
Hence from \eqref{f_psi_phi} it follows that for $k\in\left\lbrace 1,\ldots,n\right\rbrace$\\
\begin{equation*}
\psi_k = \frac{\pi k}{n+1} - \frac{2\pi k}{n+1}\sum\limits_{j=1}^n t_j + \sum\limits_{j=1}^k t_j.
\end{equation*}
But each of $t_k$ takes a finite number of values, so the set $\Lambda_0$ is also finite. Together with Lemma \ref{l_part_lim_vec} and \eqref{f_lambda_incl} this finishes the proof of Theorem \ref{t_main_ph}.  
\begin{flushright}
$\square$
\end{flushright}

\textbf{Proof of Theorem \ref{t_main}:} \eqref{f_main} directly follows from \eqref{f_main_abs_val} and \eqref{t_main_ph}.\\
\begin{flushright}
$\square$
\end{flushright}

\section{Non-constant block case.}
In this section we consider the non-constant block case. We go through the previous subsections and discuss how the argument generalizes. In order to prove Theorem \ref{t_main_non_const} Assumption \ref{assump_zm} will be again the main technical assumption which will be eliminated in the end of the proof.\\

Solution of VDE \eqref{VDE} in the non-constant block case naturally splits into $n$ components of size $N$ each: $m=\left(m^{[1]},m^{[2]},\ldots,m^{[n]}\right)$, where $m^{[k]}\in\mathbb{C}^N$ for $k\in\left\lbrace 1,\ldots,n\right\rbrace$. This partition into blocks should not be confused with the partition of the vector $m$ described in Subsection 2.1, where $m$ was split with respect to the equivalence of components. Nevertheless, in Theorem \ref{t_ord_sol_non_const} we will see that these two partitions are related to each other.\\

\subsection{Remarks to Section 2.}

All results stated in Section 2 stay unchanged. Note that we can apply Lemma \ref{l_comb} to the matrix $S$ satisfying conditions of Theorem \ref{t_main_non_const}, because $S^{j,k}_{\nu,\tau}>0$ for all $j,k\in\left\lbrace 1,\ldots,n\right\rbrace$ and $\nu,\tau\in\left\lbrace 1,\ldots,N\right\rbrace$ such that $j+k\in\left\lbrace n,n+1\right\rbrace$ and $\nu+\tau\in\left\lbrace N,N+1\right\rbrace$. Moreover, Theorem \ref{t_sol_ord} can be strengthened in the following way:\\

\begin{theorem}\label{t_ord_sol_non_const}
Assume that the matrix $S$ satisfies conditions of Theorem \ref{t_main_non_const} and Assumption \ref{assump_zm} holds. Then for every sequence $\left\lbrace z_l\right\rbrace_{l=1}^\infty$ there exists a subsequence $\left\lbrace z_{l_q}\right\rbrace_{q=1}^\infty$ such that for every $k\in\left\lbrace 1,\ldots,n\right\rbrace$ either $\|m^{[k]}\|\sim \|m^{[k+1]}\|$ with respect to the sequence $\left\lbrace z_{l_q}\right\rbrace_{q=1}^{\infty}$ or\\
\begin{equation}
\lim_{q\to\infty} {\frac {\|m^{[k]}(z_{l_q})\|}{\|m^{[k+1]}(z_{l_q})\|}} = 0.
\label{f_diff_blocks_non_const}
\end{equation}  
At the same time all components of $m$ within each block are comparable, i.e. for every $k\in\left\lbrace 1,\ldots,n\right\rbrace$ and $\nu,\tau\in\left\lbrace 1,\ldots,N\right\rbrace$ it holds that\\
\begin{equation}
\vert m^{[k]}_\nu\vert\sim\vert m^{[k]}_\tau\vert\quad \text{WRT}\quad \left\lbrace z_{l_q}\right\rbrace_{q=1}^\infty.
\label{f_same_block_non_const}
\end{equation}
\end{theorem}

\textbf{Proof of Theorem \ref{t_ord_sol_non_const}:} \eqref{f_diff_blocks_non_const} follows from Theorem \ref{t_sol_ord}, so it is sufficient to prove \eqref{f_same_block_non_const}. Without loss of generality we have $\nu\le\tau$.  Assume the opposite, i.e. that \eqref{f_same_block_non_const} does not hold. In particular, there exists a subsequence $\left\lbrace z_{l_q}\right\rbrace_{q=1}^\infty$ of the sequence $\left\lbrace z_l\right\rbrace_{l=1}^\infty$ such that\\
\begin{equation}
\lim_{q\to\infty} {\frac{\vert m^{[k]}_{\tau}(z_{l_q})\vert}{\vert m^{[k]}_\nu(z_{l_q})\vert}} = 0.
\label{f_sol_ord_non_const_1}
\end{equation}
Let $\pi\in S_{nN}$ be the transposition interchanging $\left((k-1)N+\nu\right)$ and $\left((k-1)N+\tau\right)$. Recall that the matrix $S^\pi$ is obtained from $S$ by acting with the permutation $\pi$ on the rows and  columns of $S$. Note, that $S^\pi$ still satisfies conditions of Theorem \ref{t_main_non_const} and Assumption \ref{assump_zm}. Therefore, applying Theorem \ref{t_sol_ord} to $S^\pi$ and subsequence $\left\lbrace z_{l_q}\right\rbrace_{q=1}^\infty$ we obtain that there exists a sub-subsequence $\left\lbrace z_{l_{q_u}}\right\rbrace_{u=1}^\infty$ of $\left\lbrace z_{l_q}\right\rbrace_{q=1}^\infty$ such that the sequence\\
\begin{equation}
\left\lbrace \frac{\vert m^{[k]}_{\nu}(z_{l_{q_u}})\vert}{\vert m^{[k]}_\tau(z_{l_{q_u}})\vert}\right\rbrace_{u=1}^\infty
\end{equation}
is bounded, which is contradicting to \eqref{f_sol_ord_non_const_1}. 
\begin{flushright}
$\square$
\end{flushright}

Denote the subsequence $\left\lbrace z_{l_q}\right\rbrace_{q=1}^\infty$  from Theorem \ref{t_ord_sol_non_const} by $\left\lbrace w_q\right\rbrace_{q=1}^\infty$ for short. Further in the proof of Theorem \ref{t_main_non_const} we will use this notation.\\

\subsection{Remarks to Section 3.}
Up to now the constants $c$, $C$ in Definition \ref{def_equiv} depended on the entries of $S$ (more precisely, on the lower and upper bounds of entries of $S$). Now we reinterpret notion of equivalence with respect to the sequence such that the constants $c$ and $C$ can depend on $n$ but are independent of $N$.\\ 
Further $\mathscr{C}$, $\mathscr{C}_1$ and $\mathscr{C}_2$ will be the positive constants which depend only on $\mathfrak{c}$ and $\mathfrak{C}$, but can change their value from line to line. Thus, now Definition \ref{def_equiv} takes the following form: 
\begin{defin}
For functions $f,g:\mathbb{H}\rightarrow \mathbb{C}\backslash \left\lbrace 0\right\rbrace$ and a sequence $\left\lbrace z_l\right\rbrace_{l=1}^\infty\subset \mathbb{H}$ converging to zero we will say that $f$ is equivalent to $g$ with respect to the sequence $\left\lbrace z_l\right\rbrace_{l=1}^\infty$ if there exist $\mathscr{C}_1$ and $\mathscr{C}_2$ such that\\
\begin{equation}
\mathscr{C}_1\le \frac{\vert f(z_l)\vert}{\vert g(z_l)\vert}\le \mathscr{C}_2
\end{equation}
for big enough $l$. Equivalence of functions $f$ and $g$ will be still denoted by $f\sim g$.\\
\end{defin}

All results stated in Section 3 stay unchanged. From Proposition \ref{prop_eq_phases} and Theorem \ref{t_ord_sol_non_const} it follows that\\

\begin{prop}\label{prop_eq_phases_non_const}
Assume that the matrix $S$ satisfies conditions of Theorem \ref{t_main_non_const} and Assumption \ref{assump_zm} holds. Then for every $k\in\left\lbrace 1,\ldots,n\right\rbrace$ and $\nu,\tau\in\left\lbrace 1,\ldots,N\right\rbrace$ we have\\
\begin{equation}
\lim_{q\to\infty} {\frac{\operatorname{arg}{m^{[k]}_\nu\left( w_q\right)}}{\operatorname{arg}{m^{[k]}_\tau\left( w_q\right)}}} = 0.
\end{equation}
\end{prop}

Now we will obtain from Proposition \ref{prop_eq_phases_non_const} the upper and the lower bounds for ratios of components of $m$ inside each block.

\begin{prop}\label{prop_ratio_inside_blocks}
Assume that the matrix $S$ satisfies conditions of Theorem \ref{t_main_non_const} and Assumption \ref{assump_zm} holds. Then for every $k\in\left\lbrace 1,\ldots,n\right\rbrace$ and $\nu,\tau\in\left\lbrace 1,\ldots,N\right\rbrace$ we have\\
\begin{equation}
\frac{ m^{[k]}_\nu(w_q)}{ m^{[k]}_\tau(w_q)} \sim 1 \quad \text{WRT} \quad \left\lbrace w_q\right\rbrace_{q=1}^\infty.
\label{f_ratio_non_constant}
\end{equation}
\end{prop}

\textbf{Proof of Proposition \ref{prop_ratio_inside_blocks}:} Fix $k\in\left\lbrace 1,\ldots,n\right\rbrace$ and $\nu,\tau\in\left\lbrace 1,\ldots,N\right\rbrace$. Dividing $\left((k-1)N+\tau\right)$-th equation of \eqref{VDE} by the $\left((k-1)N+\nu\right)$-th one we obtain\\
\begin{equation}
\frac{m^{[k]}_\nu(w_q)}{m^{[k]}_\tau(w_q)} = \frac{w_q+\sum\limits_{j\in\Theta_1}\sum\limits_{\xi = 1}^N S^{k,j}_{\tau,\xi}m^{[j]}_\xi(w_q) + \sum\limits_{j\in\Theta_2}\sum\limits_{\xi = 1}^N S^{k,j}_{\tau,\xi}m^{[j]}_\xi(w_q)} {w_q+\sum\limits_{j\in\Theta_1}\sum\limits_{\xi = 1}^N S^{k,j}_{\nu,\xi}m^{[j]}_\xi(w_q) + \sum\limits_{j\in\Theta_2}\sum\limits_{\xi = 1}^N S^{k,n+1-k}_{\nu,\xi}m^{[j]}_\xi(w_q)},
\label{f_ratio_non_constant_1}
\end{equation}
where\\
\begin{equation*}
\Theta_2 = \left\lbrace j\in\left\lbrace 1,\ldots,n\right\rbrace\mid \left\|m^{[j]}\right\| \sim \left\|m^{[n+1-k]}\right\|\quad \text{and}\quad S^{k,j}\neq 0\right\rbrace,
\end{equation*}
\begin{equation}
\Theta_1 = \left\lbrace 1,\ldots ,n\right\rbrace \backslash \Theta_2.
\label{def_Theta}
\end{equation}
From Assumption \ref{assump_zm} and Lemma \ref{l_summand_order_1} it follows that $\frac{z}{\|m(z)\|}$ tends to zero as $z$ tends to zero from the upper half-plane. Moreover, from Proposition \ref{prop_eq_phases} we obtain that there is no  cancellation in the sums over $\Theta_2$. Finally, from Theorem \ref{t_ord_sol_non_const} and definition \eqref{def_Theta} we have that $\sum\limits_{\Theta_1}\ll \sum\limits_{\Theta_2}$ both in numerator and denominator as $q$ tends to infinity. Therefore, it holds that\\
\begin{equation}
\left\vert\frac{m^{[k]}_\nu(w_q)}{m^{[k]}_\tau(w_q)}\right\vert = \frac{\sum\limits_{j\in\Theta_2}\sum\limits_{\xi = 1}^N S^{k,j}_{\tau,\xi}\left\vert m^{[j]}_\xi(w_q)\right\vert} {\sum\limits_{j\in\Theta_2}\sum\limits_{\xi = 1}^N S^{k,n+1-k}_{\nu,\xi}\left\vert m^{[j]}_\xi(w_q)\right\vert}(1+o(1)) \quad \text{as} \quad q\to\infty.
\label{f_ratio_non_constant_2}
\end{equation}
Replacing in \eqref{f_ratio_non_constant_2} all entries of $S$ by $\frac{\mathfrak{C}}{N}$ in numerator and by $\frac{\mathfrak{c}}{N}$ in denominator we obtain the upper bound:\\
\begin{equation}
\left\vert \frac{ m^{[k]}_\nu(w_q)}{ m^{[k]}_\tau(w_q)}\right\vert <\frac{\mathfrak{C}}{\mathfrak{c}}(1+o(1))<\mathscr{C}.
\end{equation}
The lower bound is obtained analogously.\\
\begin{flushright}
$\square$
\end{flushright}

Propositions \ref{prop_eq_phases_non_const} and \ref{prop_ratio_inside_blocks} allow us to reduce the non-constant block case to the constant block one in the following way:\\

\begin{defin}\label{def_VDE_like}
Let $\left\lbrace z_l\right\rbrace_{l=1}^\infty\subset\mathbb{H}$ be a sequence converging to zero. We will say that\\
\begin{equation}
-\frac{1}{\hat{m}(z_l)} = \omega(z_l)z_l + \hat{S}(z_l)\hat{m}(z_l)
\label{VDE_like}
\end{equation}
is a VDE-like equation with respect to the sequence $\left\lbrace z_l\right\rbrace_{l=1}^\infty$ if $\omega$, $\hat{m}$ and $\hat{S}$ are $n$-vector-functions and an $n\times n$ matrix correspondingly, which meet the following conditions:\\
\begin{enumerate}
\item{$\hat{S}(z_l)$ is a real-symmetric matrix with complex entries for all $l\in\mathbb{N}$.}\\
\item{For $j,k\in\left\lbrace 1,\ldots,n\right\rbrace$ either $\hat{s}_{j,k}(z_l)=0$ for all $l\in\mathbb{N}$ or $\hat{s}_{j,k}\sim 1$ WRT $\left\lbrace z_l\right\rbrace_{l=1}^\infty$.}\\
\item{$\operatorname{arg}{\hat{s}_{j,k}(z_l)}\to 0$ as $l\to\infty$ for all $j,k\in\left\lbrace 1,\ldots,n\right\rbrace$ such that $\hat{s}_{j,k}\sim 1$ WRT $\left\lbrace z_l\right\rbrace_{l=1}^\infty$.}\\
\item{$\omega_k(z_l)\sim 1$ WRT $\left\lbrace z_l\right\rbrace_{l=1}^\infty$ for all $k\in\left\lbrace 1,\ldots,n\right\rbrace$.}\\
\item{$\operatorname{arg}{\omega_k(z_l)}\to 0$ as $l\to\infty$ for all $k\in\left\lbrace 1,\ldots,n\right\rbrace$.}\\
\end{enumerate}
\end{defin}

\begin{theorem}\label{t_reduction}
Assume that the matrix $S$ satisfies conditions of Theorem \ref{t_main_non_const} and Assumption \ref{assump_zm} holds. Then there exists a VDE-like equation \eqref{VDE_like} with respect to the sequence $\left\lbrace w_q\right\rbrace_{q=1}^\infty$ such that $\left(m^{[k]}_1\right)_{k=1}^n$ solves \eqref{VDE_like} and for all $j,k\in\left\lbrace 1,\ldots n\right\rbrace$ we have:\\
\begin{equation}
\left(\hat{s}_{j,k}(w_q)=0 \quad \forall l\in\mathbb{N}\right) \quad \text{iff} \quad S^{j,k}=0.
\end{equation}
\end{theorem}

\textbf{Proof of Theorem \ref{t_reduction}:} 

\textbf{Step 1. Construction of VDE-like equation.} For $k\in \left\lbrace 1,\ldots,n\right\rbrace$ we sum up the equations of \eqref{VDE} corresponding to the $k$-th row of blocks of the matrix $S$:\\
\begin{equation}
-N = z\sum\limits_{\nu=1}^N m^{[k]}_\nu + \sum\limits_{j=1}^{n+k-1}\sum\limits_{\nu,\tau=1}^n S^{k,j}_{\nu,\tau} m^{[k]}_\nu m^{[j]}_\tau,
\end{equation}
that is equivalent to\\
\begin{equation}
-1 = \left(\frac{1}{N}\sum\limits_{\nu=1}^N\frac{m^{[k]}_\nu}{m^{[k]}_1}\right)zm^{[k]}_1 + \sum\limits_{j=1}^{n+k-1}\left( \sum\limits_{\nu,\tau=1}^N\frac{1}{N}S^{k,j}_{\nu,\tau}\frac{m^{[k]}_\nu}{m^{[k]}_1}\cdot\frac{m^{[j]}_\tau}{m^{[j]}_1}\right) m^{[k]}_1 m^{[j]}_1.
\end{equation}
Set\\
\begin{equation}
\omega_k(w_q):=\frac{1}{N}\sum\limits_{\nu=1}^N\frac{m^{[k]}_\nu(w_q)}{m^{[k]}_1(w_q)},\quad \forall k\in\left\lbrace 1,\ldots,n\right\rbrace
\label{def_omega}
\end{equation}
and\\
\begin{equation}
\hat{s}_{j,k}:= \sum\limits_{\nu,\tau=1}^N\frac{1}{N}S^{k,j}_{\nu,\tau}\frac{m^{[k]}_\nu(w_q)}{m^{[k]}_1(w_q)}\cdot\frac{m^{[j]}_\tau(w_q)}{m^{[j]}_1(w_q)}, \quad \forall j,k\in\left\lbrace 1,\ldots,n\right\rbrace.
\label{def_s_hat}
\end{equation}
So, $\hat{m}=\left(m^{[1]}_1,\ldots, m^{[n]}_1\right)$ solves \eqref{VDE_like} corresponding to the constructed functions $\omega$ and $S$.\\

\textbf{Step 2. Verification of properties of the constructed equation.}\\

\textbf{(1)} From the definition \eqref{def_s_hat} and the real symmetry of the matrix $S$ it follows that $\hat{S}(z_l)$ is also real symmetric for every $l\in\mathbb{N}$.\\

\textbf{(2)} If for some $j,k\in\left\lbrace 1,\ldots,n\right\rbrace$ it holds that $S^{j,k}=0$, then from \eqref{def_s_hat} it follows that $\hat{s}_{j,k}(w_q)=0$ for all $q\in\mathbb{N}$. Let us now consider $j,k\in\left\lbrace 1,\ldots,n\right\rbrace$ such that $S^{j,k}\neq 0$. From Proposition \ref{prop_eq_phases_non_const} we have that there is no cancellation in \eqref{def_s_hat}. Thus\\
\begin{equation}
\left\vert \hat{s}_{j,k}(w_q)\right\vert = \sum\limits_{\nu,\tau=1}^N\frac{1}{N}S^{k,j}_{\nu,\tau}\left\vert\frac{m^{[k]}_\nu(w_q)}{m^{[k]}_1(w_q)}\right\vert\cdot\left\vert\frac{m^{[j]}_\tau(w_q)}{m^{[j]}_1(w_q)}\right\vert (1+o(1)) \quad \text{as} \quad q\to\infty.
\end{equation}
Therefore, from Proposition \ref{prop_ratio_inside_blocks} we have\\
\begin{equation}
\left\vert \hat{s}_{j,k}(w_q)\right\vert \le \sum\limits_{\nu,\tau=1}^N\frac{1}{N} \frac{\mathfrak{C}}{N}\mathscr{C}\le \mathscr{C}_1.
\end{equation}
Similarly we obtain that $\left\vert \hat{s}_{j,k}(w_q)\right\vert\ge\mathscr{C}_2$. So, $\hat{s}_{j,k}(w_q)\sim 1$ WRT $\left\lbrace w_q\right\rbrace_{q=1}^\infty$.\\

\textbf{(3)} From Proposition \ref{prop_eq_phases_non_const} it follows that the argument of each summand in \eqref{def_s_hat} tends to zero as $q$ tends to infinity. Therefore, $\operatorname{arg}{\hat{s}(w_q)}\to 0$ as $l\to\infty$ for every non-zero block $S^{j,k}$.\\
Verification of the similar properties of $\omega(z_l)$ is exactly the same as in (2) and (3) above.\\
\begin{flushright}
$\square$
\end{flushright}

Further in the proof of Theorem \ref{t_main_non_const} we consider \eqref{VDE_like} instead of \eqref{VDE}, where \eqref{VDE_like} is obtained from \eqref{VDE} .\\

\subsection{Remarks to Section 4 and the end of the proof of Theorem \ref{t_main_non_const}.}
Theorem \ref{t_main_seq} consists of two parts: equivalence of the solutions of \eqref{VDE} to the certain powers of $z$ and existence of the limit constants $\left\lbrace c_k\right\rbrace_{k=1}^n$. Note that in Section 4 this equivalence was obtained in the old sense and in the current subsection we will prove that it still holds in the new sense of equivalence. The first part literally generalizes to the non-constant block case, but the second part cannot be generalized directly, because the entries of the matrix $\hat{S}$ are not constants. However, the same argument will give the estimates for the partial limits of the functions
\begin{equation}
\vert m_k(z)\vert\cdot\vert z\vert ^{-\left(1-\frac{2k}{n+1}\right)}
\end{equation}
as $z$ tends to zero from the complex upper half-plane.

\begin{theorem}\label{t_main_seq_non_const}
Assume that the matrix $S$ meets conditions of Theorem \ref{t_main_non_const} and Assumption \ref{assump_zm} holds. Then with the notations of Theorem \ref{t_ord_sol_non_const} for every $k\in\left\lbrace 1,\ldots,n\right\rbrace$ we have
\begin{equation}
\vert m_k(w_q)\vert\cdot\vert w_q\vert ^{-\left(1-\frac{2k}{n+1}\right)} \sim 1 \quad \text{WRT}\quad \left\lbrace w_q\right\rbrace_{q=1}^\infty.
\end{equation}
\end{theorem}

\textbf{Proof of Theorem \ref{t_main_seq_non_const}:} As it was noted above, for some positive constants $c$ and $C$ we have\\
\begin{equation}
c<\vert m_k(w_q)\vert\cdot\vert w_q\vert ^{-\left(1-\frac{2k}{n+1}\right)}<C
\end{equation}
for big enough $q$, but we do not know how $c$ and $C$ depend on $N$. In order to obtain the estimates for these constants we repeat the argument from the proof of Theorem \ref{t_main_seq_non_const}: consider an arbitrary subsequence $\left\lbrace w_{q_u}\right\rbrace_{u=1}^\infty$ of the sequence $\left\lbrace w_q\right\rbrace_{q=1}^\infty$ such that the following limits exist:\\
\begin{equation}
\hat{c}_k:=\lim_{u\to\infty} {\left(\vert m_k(w_{q_u}\vert\cdot\vert w_{q_u}\vert ^{-\left(1-\frac{2k}{n+1}\right)}\right)}.
\end{equation}
Up to extraction of a subsequence we can assume that there also exist the limits\\
\begin{equation}
\hat{s}^0_{j,k}:=\lim_{u\to\infty} \hat{s}_{j,k}(w_{q_u}),\quad \forall j,k\in\left\lbrace 1,\ldots,n\right\rbrace
\end{equation}
and\\
\begin{equation}
\omega^0_k:=\lim_{u\to\infty} \omega_k(w_{q_u}),\quad \forall k\in\left\lbrace 1,\ldots,n\right\rbrace.
\end{equation}
From the definition of VDE-like equation it follows that $\hat{s}^0_{j,k}$ and $\omega^0_k$ are positive for all $j,k\in\left\lbrace 1,\ldots,n\right\rbrace$. So, as in the proof of Theorem \ref{t_main_seq} we can show that $\left\lbrace \log{\hat{c}_k} \right\rbrace_{k=1}^n$ solves the linear system \eqref{f_logc_sys}, i.e.\\
\begin{equation}
\begin{cases}
\log{\hat{c}_k}+\log{\hat{c}_{n+1-k}}=-\log{s^0_{k,n+1-k}}, &1\le k\le \left\lceil\frac{n}{2}\right\rceil\\
\log{\hat{c}_1}+\log{\hat{c}_{n-1}} - \log{\hat{c}_n} = \log{\omega^0_n}-\log{s^0_{1,n-1}}\\
\log{\hat{c}_k}+ \log{\hat{c}_{n-k}} - \log{\hat{c}_{n+1-k}} - \log{\hat{c}_{k-1}}=\log{s^0_{n+1-k,k-1}}-\log{s^0_{k,n-k}}, &2\le k\le \left[\frac{n}{2}\right].\\
\end{cases}
\label{f_logc_sys_non_const}
\end{equation}
Rewrite the system \eqref{f_logc_sys_non_const} in the following way: $A\log{\hat{c}}=b$, where $\log{\hat{c}}$ is a vector whose components are $\log{\hat{c}_k}$ for $k\in\left\lbrace 1,\ldots,n\right\rbrace$. In the proof of Theorem \ref{t_main_seq} it was shown that the matrix $A$ is invertible. Note that $A$ does not depend on $N$ and for vector $b$ we have an upper bound which does not depend on $N$. Therefore, $\|\log{\hat{c}}\|\le \|A^{-1}\|\cdot \|b\|\le\mathscr{C}$, so $\hat{c}_k\sim 1$ for every $k\in\left\lbrace 1,\ldots,n\right\rbrace$. 
\begin{flushright}
$\square$
\end{flushright}

Results stated in Sections 5 and 6 can be literally generalized to the case of VDE-like equation. Therefore, Theorem \ref{t_main_non_const} is proved.\\


\begin{thebibliography}{5}
\bibitem{Erd}
L. Erd{\H o}s, {\it The matrix Dyson equation and its applications for random matrices.}
{\it in} Random matrices, 75--158, IAS/Park City Math. Ser., 26, Amer. Math. Soc., Providence, RI, 2019.
\bibitem{Frob}
Roger A. Horn and Charles R. Johnson. Matrix Analysis. 2nd. Cambridge University Press, 2012
\bibitem{Asb}
Asbjorn Baekgaard Lauristen. The Vector Dyson equation for a class of\\
skew-triangular block matrices. 2021. Unpublished manuscript see\\ 
{\footnotesize $https://www2.ist.ac.at/fileadmin/download/ReportAsbjorn\_Baekgaard\_LAURITSEN\_.pdf $}
\bibitem{Erd1}
O. Ajanki, L. Erd{\H o}s,  T. Kr\"uger, {\it Quadratic vector equations on complex upper half-plane.}
Memoirs of Amer. Math. Soc. vol. {\bf 261} No. 1261 (2019)
\bibitem{Torb}
T.Kr\"uger, D.Renfrew, Singularity degree of structured random matrices, preprint on arXiv.\\
\end{thebibliography}
\end{document}